\documentclass[11pt, oneside]{article}
\usepackage[utf8]{inputenc}
\usepackage{amsmath, amsthm, amssymb}
\usepackage[dvipsnames]{xcolor}
\usepackage{enumerate, comment, multirow, parskip}
\usepackage{dsfont}
\usepackage{enumitem}
\usepackage{mathtools}

\usepackage{subcaption}
\usepackage{todonotes}
\usepackage{multicol}
\usepackage{tikz}
\usetikzlibrary{decorations.markings}
\usetikzlibrary{decorations.pathmorphing}
\usetikzlibrary{patterns}
\usepackage{standalone}
\usepackage{mathtools,mdframed}
\usepackage{adjustbox}
\usepackage[left=1in,right=1in,top=1in,bottom=1in,bindingoffset=0cm]{geometry}

\usepackage[normalem]{ulem}

\usepackage{hyperref}
\hypersetup{colorlinks=true,linkcolor=blue,anchorcolor=blue,citecolor=green,filecolor=blue,urlcolor=blue,bookmarksnumbered=true,pdfview=FitB}

\usepackage[capitalise]{cleveref}

\definecolor{green}{RGB}{50,205,50}
\definecolor{yellow}{RGB}{240,240,10}

\newtheorem{theorem}{Theorem}[section]
\newtheorem{lemma}[theorem]{Lemma}
\newtheorem{corollary}[theorem]{Corollary}

\newtheorem{conjecture}[theorem]{Conjecture}
\newtheorem{obs}{Observation}

\newtheorem*{theorem*}{Theorem}
\newtheorem*{corollary*}{Corollary}

\theoremstyle{definition}

\newtheorem{example}{Example}

\theoremstyle{definition}

\newtheorem{claim}{Claim}
\newenvironment{proofc}{\begin{proof}[Proof of Claim]}{\end{proof}}

\makeatletter
\pgfdeclarepatternformonly{fine north east lines}
  {\pgfqpoint{.5pt}{.5pt}}
  {\pgfqpoint{2pt}{2pt}} 
  {\pgfqpoint{2pt}{2pt}} 
  {
    \pgfsetlinewidth{.2pt}
    \pgfpathmoveto{\pgfqpoint{0pt}{0pt}}
    \pgfpathlineto{\pgfqpoint{4pt}{4pt}}
    \pgfusepath{stroke}
  }
\makeatother

\tikzset{
  finer stripes/.style={
    pattern=fine north east lines,
    pattern color=black
  }
}

\title{Extending total colorings in planar graphs}

\author{
Owen Henderschedt\footnote{Auburn University, Department of Mathematics and Statistics, Auburn U.S.A.
  Email: {\tt olh0011@auburn.edu}.}
\and
Jessica McDonald\footnote{Auburn University, Department of Mathematics and Statistics, Auburn U.S.A.
  Email: {\tt mcdonald@auburn.edu}.   
	Supported in part by Simons Foundation Grant \#845698 and NSF grant DMS-2452103.}}
\date{}

\begin{document}
\maketitle

\begin{abstract}
We initiate the study of total-coloring extensions, and focus our attention on planar graphs, asking: ``When can a total-$k$-coloring of some subgraph $H$ of a planar graph $G$ be extended to a total-$k$-coloring of $G$?'' 
We prove that if  $H$ is a matching, then any total-$(\Delta+3)$-coloring of $H$ in $G$ extends to $G$ provided $\Delta\geq 28$; this number of colors is best-possible without introducing a distance condition on $H$. We also prove that if $H$ is a set of distance-3 cliques then any total-$(\Delta+1)$-coloring of $H$ extends to $G$ provided  $\Delta\geq 27$; this distance condition cannot be lowered. 
\end{abstract}

\section{Introduction}

Once a coloring result is known, it is natural to ask how robust that result is via a \textit{coloring extension problem}. Namely, if we know that a graph $G$ is $k$-colorable we may ask: \textit{when can a $k$-coloring of some subgraph $H$ of $G$ be extended to a $k$-coloring of $G$?} Extension problems have been studied for both vertex-colorings and edge-colorings, particularly when $G$ is planar. In this paper we initiate the study of total-coloring extension problems. We prove three main results, all of which are for planar graphs, and conjecture a number of generalizations. Our results concern the cases when $H$ is a matching, when $H$ is more generally a disjoint set of cliques, and when $H$ has bounded maximum degree. Each of these results is motivated by work that has been done on vertex-coloring and edge-coloring extension problems, and as well as total and list-total coloring results. 

In this paper all graphs are simple; we refer the reader to \cite{west} for terms not defined here. Given a graph $G$, a \emph{total-$k$-coloring} is an assignment of the colors $\{1, 2, \ldots, k\}=[k]$ to the vertices and edges of $G$ such that the vertices induce a $k$-coloring, the edges induce a $k$-edge-coloring, and no edge can receive the same color as either of its endpoints. The \emph{total chromatic number} of a graph $G$, denoted  $\chi_T(G)$, is the smallest $k$ such that there exists a total-$k$-coloring of $G$. By considering the edges incident with a vertex of maximum degree $\Delta(G)$, along with the vertex itself, we see that $\chi_T(G)\geq \Delta(G)+1$. In the 1960's, Behzad \cite{Behzad-total-coloring} and Vizing \cite{Vizing-total-coloring} independently proposed the \emph{Total Coloring Conjecture} which says that $\chi_T(G)$ is at most one more than this lower bound.

\begin{conjecture}[The Total Coloring Conjecture]\label{conj:TCC}
For any graph $G$,  $\chi_T(G)\leq \Delta(G)+2$.
\end{conjecture}

Conjecture \ref{conj:TCC} has been verified for all planar graphs $G$ with $\Delta(G)\neq 6$ through a series of papers (see \cite{KOSTOCHKA-MAX-DEG-5,SANDERS-DEG-7}). In fact, it is known that $\chi_T(G)= \Delta(G)+1$ for planar graphs with $\Delta(G)\geq 9$ (see \cite{TYPE-1-MAX-DEG-9}), although this is false when $\Delta(G)\leq 3$.  For much more on the Total Coloring Conjecture, both in the planar case and in general, the interested reader is referred to 
\cite{GNS} and \cite{JT}. 

Suppose we have a graph $G$ with a subgraph $H$ and $\varphi$ is a total-$k$-coloring of $H$. We want to know whether $\varphi$ can be \emph{extended} to $G$, that is, whether there exists a total-$k$-coloring of $G$ which agrees with $\varphi$ on all the vertices and edges of $H$. Of course, if there are vertices $u$ and $v$ which are adjacent in $G$ but non-adjacent in $H$ and moreover $\varphi(u)=\varphi(v)$, then we cannot hope to extend $\varphi$ to $G$. We say that $\varphi$ is a \emph{total-$k$-coloring of $H$ in $G$} if no such vertices $u, v$ exist, that is, if the coloring $\varphi$ is proper also with respect to $G$. In general, if $\varphi$ is a total-$k$-coloring of $H$ in $G$, then any uncolored vertex in $G$ sees at most $2\Delta(G)$ colors from $\varphi$ (all adjacent vertices and incident edges) and any uncolored edge in $G$ sees at most $2\Delta(G)$ colors ($2\Delta(G)-2$ adjacent edges and $2$ endpoints). So we get the following basic observation by coloring greedily.

\begin{obs}\label{obs: greedy} If $G$ is a graph with maximum degree $\Delta$ and $H$ is a subgraph of $G$, then any total-$(2\Delta+1)$-coloring of $H$ in $G$ can be extended to $G$.
\end{obs}

Without putting any restrictions on $H$, the $(2\Delta+1)$ in Observation \ref{obs: greedy} is best possible due to \cref{ex: greedy} (given in Section 2).  So we must make restrictions on the precolored subgraph $H$ if we want to lower the number of colors beyond the $(2\Delta+1)$, and certainly if we want to approach the $(\Delta+2)$ suggested by the Total Coloring Conjecture. Our starting point for this is matchings. While a matching is typically defined as a set of edges within a graph whose endpoints are vertex disjoint, we extend this notion here and say that a subgraph $H$ of a graph $G$ is a \emph{matching} if it induces a set of matching edges in $H$, i.e., if $H$ consists of vertex-disjoint copies of $K_2$ within $G$.

\begin{conjecture}\label{conj: ALL-G-Delta+3}
Let $G$ be a graph with maximum degree $\Delta$ and let $H$ be a subgraph of $G$. If $H$ is a matching, then any total-$(\Delta+3)$-coloring of $H$ in $G$ can be extended to $G$. 
\end{conjecture}

\begin{theorem}\label{thm: D+3-no-distance}
Let $G$ be a planar graph with maximum degree $\Delta\geq 28$ and let $H$ be a subgraph of $G$. If $H$ is a matching, then any total-$(\Delta+3)$-coloring of $H$ in $G$ can be extended to $G$. 
\end{theorem}

The $(\Delta+3)$ in \cref{conj: ALL-G-Delta+3} and \cref{thm: D+3-no-distance} (the first main result of this paper), is best possible due to \cref{ex: sharp distance2} (given in Section 2). In fact, this example shows that  $(\Delta+3)$ cannot even be replaced by $(\chi_T(G)+1)$.

Given a graph $G$ with a subgraph $H$ that is a matching, we say that $H$ is 
\emph{distance-$\ell$} if the distance between any two vertices in distinct components of $H$ is at least $\ell$ in $G$. The matching in \cref{ex: sharp distance2} is distance-2, and we conjecture that this distance is the barrier for extending from totally-$(\Delta+2)$-colored matchings. 

\begin{conjecture}\label{conj: Delta+2-strong}
Let $G$ be a graph with maximum degree $\Delta$ and let $H$ be a subgraph of $G$. If $H$ is a distance-3 matching, then any total-$(\Delta+2)$-coloring of $H$ can be extended to $G$. 
\end{conjecture}

While we can't hope for fewer than $\Delta+2$ colors in general, we can actually do better than Conjecture \ref{conj: Delta+2-strong} in the planar case when the maximum degree is large enough.

\begin{theorem}\label{thm: Delta+1-match}
Let $G$ be a planar graph with maximum degree $\Delta\geq 27$ and let $H$ be a subgraph of $G$. If $H$ is a distance-3 matching, then any total-$(\Delta+1)$-coloring of $H$ can be extended to $G$. 
\end{theorem}

Note that we do not need to say that the coloring of $H$ is `in $G$' in Conjecture \ref{conj: Delta+2-strong} or Theorem Theorem \ref{thm: Delta+1-match}, since this is implicit from the distance condition. We will actually get Theorem \ref{thm: Delta+1-match} as a corollary of our third main result, which we will discuss towards the end of this section. In general, Conjecture \ref{conj: Delta+2-strong} and Theorem \ref{thm: Delta+1-match} are total-coloring analogs of the following conjecture and theorem of Edwards, Gir\~ao, van den Heuvel, Kang, Puleo, and Sereni \cite{EDWARDS}.
             
\begin{conjecture}[Edwards, et. al. \cite{EDWARDS}]\label{conj: Edwards-edge-extension}
Let $G$ be a graph with maximum degree $\Delta$, and let $H$ be a subgraph of $G$. If $H$ is a distance-2 matching, then any $(\Delta+1)$-edge-coloring of $H$ can be extended to $G$. 
\end{conjecture} 

\begin{theorem}[Edwards, et. al. \cite{EDWARDS}]\label{thm: Edwards}
Let $G$ be a planar graph with maximum degree $\Delta\geq 17$, and let $H$ be a subgraph of $G$. If $H$ is a distance-2 matching, then any $(\Delta+1)$-edge-coloring of $H$ can be extended to $G$. 
\end{theorem} 

To be clear, by `extending the $(\Delta+1)$-edge-coloring $\varphi$ of $H$ to $G$', we mean that $G$ has a $(\Delta+1)$-edge-coloring which agrees with $\varphi$ on the edges of $H$. Conjecture \ref{conj: Edwards-edge-extension} is a strengthening of an older conjecture of Albertson and Moore in \cite{ALBERTSON}, and there are other results supporting it by 
Gir\~ao and Kang \cite{GIRAO} and Cao, Chen, Qi, and Shan \cite{CAO} (including for a multigraph version). In comparing Conjectures \ref{conj: Delta+2-strong} and \ref{conj: Edwards-edge-extension}, observe that while Vizing's Theorem \cite{VIZING} guarantees that all graphs have a $(\Delta+1)$-edge-coloring, Conjecture \ref{conj: Delta+2-strong} would be stronger than the Total Coloring Conjecture.

The second main result of this paper requires no particular structure for the subgrapah $H$, but rather requires $H$ to have bounded maximum degree. 

\begin{theorem}\label{thm:Planar-Delta+t}
Let $G$ be a planar graph with maximum degree $\Delta$ and let $H$ be a subgraph of $G$. If $H$ has maximum degree $d$, then any $(\Delta+d+6)$-total-coloring of $H$ in $G$ can be extended to $G$. 
\end{theorem}

Theorem \ref{thm:Planar-Delta+t} is a total-coloring analog of the following result.

\begin{theorem}[Harrelson, McDonald, Puleo \cite{HMP}]\label{thm:HMP}
Let $G$ be a planar graph with maximum degree $\Delta$ and let $H$ be a subgraph of $G$. If $H$ has maximum degree $d$, then any $(\Delta+d+4)$-edge-coloring of $H$ can be extended to $G$. 
\end{theorem}

Example \ref{ex: greedy} shows that the $(\Delta+d+6)$ in Theorem \ref{thm:Planar-Delta+t} cannot be lowered beyond $(\Delta+d+2)$; in that example $d=\Delta-1$ and a total-$(\Delta+d+1)$-coloring of $H$ in $G$ cannot be extended to $G$. Harrelson, McDonald, and Puleo have a similar example in \cite{HMP} which shows that $(\Delta+d)$ would be best-possible for Theorem \ref{thm:HMP}; in fact they prove this strongest result whenever $\Delta\geq 16+d$. By making a similar assumption on $\Delta$, it may be possible to get a total-coloring analog with $(\Delta+d+2)$, but we leave this to future work. In fact, our main purpose in including Theorem \ref{thm:Planar-Delta+t} in this paper is to apply it with $d=1$ as a base case in our proof for Theorem \ref{thm: D+3-no-distance} (there, the value $(\Delta+7)$ is reduced to $(\Delta+3)$ at the cost of requiring $\Delta$ to be sufficiently large).

While the literature on extending edge-colorings focuses on extending from (sufficiently distanced) matchings; the literature on extending vertex-colorings is rich with results for extending from (sufficiently distanced) cliques. A variety of such results may be found in \cite{PAINT-CORNER, BROOKS, ALBERTSON, HUTCHINSON-MOORE, OJIMA}; for our purposes we will highlight just one: Albertson \cite{PAINT-CORNER} proved that in any planar graph $G$, if a set of distance-3 vertices are 6-colored, then this can always be extended to $G$. To be clear, given a subgraph $H$ of a graph $G$, we say that $H$ is a \emph{set of distance-$\ell$ cliques} if $H$ is a set of vertex-disjoint cliques (which may be of varying sizes) and the distance between any two vertices in distinct components of $H$ is at least $\ell$. In order to give more context for our third and final main result of this paper, which is about extending total-colorings from a set of distance-3 cliques, we need to first discuss total-choosability.

All total-coloring extension problems can be rephrased as list-total-coloring problems. We say that a graph $G$ is \emph{totally-$L$-choosable} if $G$ can be totally colored so that each $x\in V(G)\cup E(G)$ is assigned a color from a list $L(x)$. If $G$ is totally-$L$-choosable for all $L$ where $|L(x)|=k$ for each $x\in V(G)\cup E(G)$, then we say that $G$ is \emph{totally-$k$-choosable}. Suppose we have a graph $G$ and subgraph $H$, and we want to extend a total-$k$-coloring $\varphi$ of $H$ in $G$. For each $x\in V(G)\cup E(G)\setminus (V(H)\cup E(H))$ we define $L(x)$ as the subset of $[k]$ remaining after removing any color incident or adjacent to $x$ under $\varphi$. If we also set $L(x)=\{\varphi(x)\}$ for all $x\in V(H)\cup E(H)$, then $G$ is totally-$L$-choosable if and only if $\varphi$ can be extended to $G$. Moreover, for \emph{any} way we extend $L$ by assigning lists $L(x)$ for $x\in V(H)\cup E(H)$, if we get that $G$ is totally-$L$-choosable, then it implies that $\varphi$ can be extended to $G$. This is because the color assigned to some $x\in V(H)\cup E(H)$ can be removed and replaced with $\varphi(x)$ without causing any conflicts.

Every planar graph $G$ with maximum degree $\Delta\geq 12$ is known to be total-$(\Delta+1)$-choosable by Borodin, Kostochka, and Woodall \cite{DELTA-PLUS-1}. Suppose that $H$ is a set of disjoint cliques in $G$ which has been totally colored by $\varphi$. If $H$ is a distance-3 matching, then any vertex or edge in $G$ sees at most two different colors of $\varphi$. If $H$ is more generally a set of distance-3 cliques, then any vertex or edge in $G$ sees at most $\omega(H)$ colors (and $\omega(G)\leq 4$ since $G$ is planar). So, using the approach described in the above paragraph, Borodin and Kostochka's total choosability result implies the following.

\begin{theorem}[Borodin, Kostochka, and Woodall \cite{DELTA-PLUS-1}]\label{PlanarDelta+3}
Let $G$ be a planar graph with maximum degree $\Delta\geq 12$ and let $H$ be a subgraph of $G$.
\begin{enumerate}
\item If $H$ is a distance-3 matching, then any total-$(\Delta+3)$-coloring of $H$ can be extended to $G$.
\item If $H$ is a set of distance-3 cliques, then any total-$(\Delta+5)$-coloring of $H$ can be extended to $G$.
\end{enumerate}
\end{theorem}

We can compare Theorem \ref{PlanarDelta+3}(1) to Conjecture \ref{conj: ALL-G-Delta+3} and Theorem \ref{thm: D+3-no-distance} where we have the same number of colors but no distance condition; we can also compare it to Conjecture \ref{conj: Delta+2-strong} and Theorem \ref{thm: Delta+1-match} where we have the same distance condition but we are using one or even two fewer colors, respectively.

Our third and final main result of this paper lowers the number of colors in Theorem \ref{PlanarDelta+3}(2) from $(\Delta+5)$ all the way down to $(\Delta+1)$, as follows.

\begin{theorem}\label{thm: D+t - all cliques}
Let $G$ be a planar graph with maximum degree $\Delta\geq 27$, and let $H$ be a subgraph of $G$. If $H$ is a set of distance-3 cliques, then any $(\Delta+1)$-total-coloring of $H$ can be extended to $G$.
\end{theorem}

We now see that Theorem \ref{thm: D+t - all cliques} immediately implies Theorem \ref{thm: Delta+1-match}. Example \ref{ex: sharp distance2} again shows that the distance-3 condition in Theorem \ref{thm: D+t - all cliques} cannot be lowered to distance-2. The $\Delta\geq 27$ condition in Theorem \ref{thm: D+t - all cliques} is surely not sharp, but it cannot be eliminated completely. In particular, not all planar graphs with maximum degree $\Delta\leq 3$ are total-$(\Delta+1)$-colorable. In case one cares a great deal about the bound on $\Delta$, we have actually proved slightly stronger versions of Theorem \ref{thm: D+3-no-distance} and Theorem \ref{thm: D+t - all cliques} than we stated in this introduction -- in particular, we show that the bound on $\Delta$ can be improved if more colors are used. Actually, we also allow precolored vertices to be present for Theorem \ref{thm: D+3-no-distance},  at no cost.

The current paper now proceeds as follows. Examples \ref{ex: greedy}, \ref{ex: sharp distance2}, and \ref{ex: distance-4} appear in Section \ref{sec:examples}. Then we prove  \cref{thm: D+t - all cliques} in Section \ref{sec: thm1.8}, Theorem \ref{thm:Planar-Delta+t} in Section \ref{sec: thm1.6}, and Theorem \ref{thm: D+3-no-distance} in Section \ref{sec: thm1.3}. 
In all three of our main proofs, when we consider an uncolored edge or vertex, we will be counting how many colors are available for it. Such counts are independent of the color palette (eg. $\{1, 2, \ldots, k\}$ for a total-$k$-coloring), and rather just require that each structure starts with a list of $k$ distinct colors. In this way, all our proofs remain valid in the more general setting of \emph{total list-coloring extensions}, in which each vertex and edge is assigned its own list of available colors and a subgraph is precolored from those lists, with the goal of extending the precoloring to the entire graph. In this paper we have chosen to frame all our results as total coloring extensions, rather than total list-coloring extensions, in order to emphasize readability and present our theorems in what we feel is their most natural setting.

\section{Examples}\label{sec:examples}

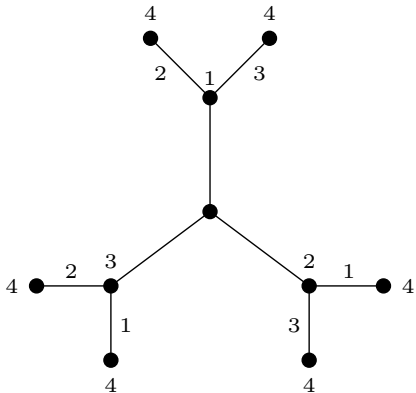
\begin{figure}[h!]
    \centering
    \begin{tikzpicture}











\begin{scope}[shift = {(5,0)}]

    \draw[] (0,0) to (0,1.1);
    \draw[] (0,0) to (-1,-.75);
    \draw[] (0,0) to (1,-.75);

    \draw[] (0,1.15) to (-.6,1.75);
    \draw[] (0,1.15) to (.6,1.75);

    \draw[] (-1,-.75) to (-1.75,-.75);
    \draw[] (-1,-.75) to (-1,-1.5);

    \draw[] (1,-.75) to (1.75,-.75);
    \draw[] (1,-.75) to (1,-1.5);

    \draw[fill = black] (0,0) circle (2pt);
    
    \draw[fill = black] (-1,-.75) circle (2pt);
    \draw[fill = black] (1,-.75) circle (2pt);
    \draw[fill = black] (0,1.15) circle (2pt);
    
    \draw[fill = black] (-1,-1.5) circle (2pt);
    \draw[fill = black] (1,-1.5) circle (2pt);
    
    \draw[fill = black] (-.6,1.75) circle (2pt);
    \draw[fill = black] (.6,1.75) circle (2pt);

    \draw[fill = black] (-1.75,-.75) circle (2pt);
    \draw[fill = black] (1.75,-.75) circle (2pt);

\node at (0,1.35) {\tiny{$1$}};
\node at (1,-.5) {\tiny{$2$}};
\node at (-1,-.5) {\tiny{$3$}};

\node at (-.5,1.4) {\tiny{$2$}};
\node at (.5,1.4) {\tiny{$3$}};

\node at (1.4,-.6) {\tiny{$1$}};
\node at (.85,-1.15) {\tiny{$3$}};

\node at (-1.4,-.6) {\tiny{$2$}};
\node at (-.85,-1.15) {\tiny{$1$}};

\node at (-.6,2) {\tiny{$4$}};
\node at (.6,2) {\tiny{$4$}};
\node at (-2,-.75) {\tiny{$4$}};
\node at (2,-.75) {\tiny{$4$}};
\node at (-1,-1.75) {\tiny{$4$}};
\node at (1,-1.75) {\tiny{$4$}};

\end{scope}

\end{tikzpicture}
    \caption{\cref{ex: greedy} with $k=3$ showing that the color palette $[2\Delta+1]$ cannot be reduced in size for guaranteeing and extension of a total-precoloring to total-colorings.}
    \label{fig:greedyexample}
\end{figure}

\begin{example}\label{ex: greedy}
Let $G$ be the tree which is rooted at vertex $v$ with $k$ children $u_1,...,u_k$, each themselves with $k-1$ children. Let $W_i=w^{(i)}_1,...w_{k-1}^{(i)}$ be the children of $u_i$. We have $\Delta(G)=k$. We now define a total-precoloring of a subgraph of $G$ as follow. Let $\varphi(u_i) =i$ for all $i\in [k]$, let $\varphi(w_j^{(i)})=k+1$ for all $i\in [k]$ and $j\in [k-1]$ and let $\varphi(u_iw_j^{(i)}) = i+j \pmod{k}$ (where $0 = k$). Note that $\varphi$ is a total-$k$-coloring of $k$ disjoint $K_{1,(k-1)}$'s. Viewing $\varphi$ as a total-$2k$-coloring (with colors $\{k+1, \ldots, 2k\}$ unused) we claim that it does not  extend to $G$. To see this, notice that every edge $vu_i$ and $v$ sees every color in $[k]$ and since $vu_i$ and $vu_j$ are pairwise adjacent for all $i,j\in [k]$, all $k+1$ items must receive a distinct color (see \cref{fig:greedyexample}).
\end{example}

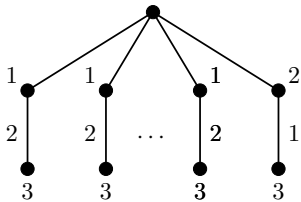
\begin{figure}[h!]
    \centering
    \begin{tikzpicture}

\draw[thick] (0,.25) to (-.75,-1); 
\draw[thick] (0,.25) to (-2,-1); 
\draw[thick] (0,.25) to (.75,-1); 
\draw[thick] (0,.25) to (2,-1); 
\draw[thick] (-.75,-2.25) to (-.75,-1); 
\draw[thick] (-2,-2.25) to (-2,-1); 
\draw[thick] (.75,-2.25) to (.75,-1); 
\draw[thick] (2,-2.25) to (2,-1);

\draw[fill = black] (0,.25) circle (3pt);
\draw[fill = black] (-.75,-1) circle (3pt);
\draw[fill = black] (-2,-1) circle (3pt);
\draw[fill = black] (.75,-1) circle (3pt);
\draw[fill = black] (2,-1) circle (3pt);
\draw[fill = black] (-.75,-2.25) circle (3pt);
\draw[fill = black] (-2,-2.25) circle (3pt);
\draw[fill = black] (.75,-2.25) circle (3pt);
\draw[fill = black] (2,-2.25) circle (3pt);

\node at (0,-1.75) {$\cdots$};
\node at (-2-.25, -1+.25) {$1$};
\node at (-2-.25, -1.675) {$2$};
\node at (-2, -2.6) {$3$};

\node at (-.75-.25, -1+.25) {$1$};
\node at (-.75-.25, -1.675) {$2$};
\node at (-.75, -2-.6) {$3$};

\node at (.75+.25, -1+.25) {$1$};
\node at (.75+.25, -1.675) {$2$};
\node at (.75, -2-.6) {$3$};

\node at (.75+.25, -1+.25) {$1$};
\node at (.75+.25, -1.675) {$2$};
\node at (.75, -2-.6) {$3$};

\node at (2+.25, -1+.25) {$2$};
\node at (2+.25, -1.675) {$1$};
\node at (2, -2-.6) {$3$};

\end{tikzpicture}
    \caption{\cref{ex: sharp distance2} with $t=4$ which shows a totally precolored distance-$2$ matching with color palette $[\Delta+2]$ which cannot extend to a total-$(\Delta+2)$-coloring}
    \label{fig:distance2-sharp}
\end{figure}

\begin{example}\label{ex: sharp distance2}
Consider the graph $G$ obtained from a star $K_{1,t}$ where each edge is subdivided once. Let $v$ be the vertex of degree $t$ with neighbors $x_1,...,x_t$ of degree $2$. Let $y_1,...,y_t$ be the leaf neighbors of $x_1,...,x_t$ respectively. Note that $\Delta(G)=t$. Now consider the total-$(t+2)$-coloring $\varphi$ of $H = \{x_iy_i: i\in [t]\}$ where $\varphi(x_i) = 1$ for all $i\in [t-1]$ and, $\varphi(x_t)=2$, $\varphi(x_iy_i)=2$ for all $i\in [t-1]$ and $\varphi(x_ty_t)=1$. Finally let $\varphi(y_i)=3$ for all $i\in [t]$ (note that the coloring of the $y_i$'s can be colored arbitrarily). See \cref{fig:distance2-sharp} for the setup of $\varphi$ and $G$. Now notice there are $t+1$ remaining items ($t$ edges and $1$ vertex) to color and they all must be colored with a distinct color. But since $1$ and $2$ are forbidden for both $v$ and $vx_i$ for all $i\in [t]$, we cannot extend $\varphi$ to a total-$(t+2)$-coloring of $G$. 
\end{example}

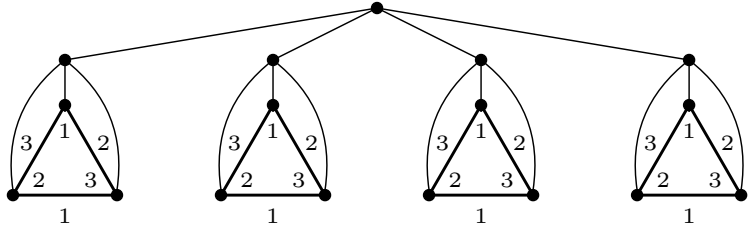
\begin{figure}[h!]
    \centering
    \begin{tikzpicture}
\begin{scope}[shift = {(0,0)}]

\draw (0,0) -- (-3,-.5);
\draw (0,0) -- (-1,-.5);
\draw (0,0) -- (1,-.5);
\draw (0,0) -- (3,-.5);

\draw[bend left = -30] (-3,-.5) to (-3.5,-1.75);
\draw[bend left = 30] (-3,-.5) to (-2.5,-1.75);
\draw[bend left = 0] (-3,-.5) to (-3,-1);

\draw[bend left = -30] (-1,-.5) to (-1.5,-1.75);
\draw[bend left = 30] (-1,-.5) to (-.5,-1.75);
\draw[bend left = 0] (-1,-.5) to (-1,-1);

\draw[bend left = -30] (1,-.5) to (.5,-1.75);
\draw[bend left = 30] (1,-.5) to (1.5,-1.75);
\draw[bend left = 0] (1,-.5) to (1,-1);

\draw[bend left = 30] (3,-.5) to (3.5,-1.75);
\draw[bend left = -30] (3,-.5) to (2.5,-1.75);
\draw[bend left = 0] (3,-.5) to (3,-1);

\draw[fill = black] (0,0) circle (1.5pt);
\draw[fill = black] (-3,-.5) circle (1.5pt);
\draw[fill = black] (-1,-.5) circle (1.5pt);
\draw[fill = black] (1,-.5) circle (1.5pt);
\draw[fill = black] (3,-.5) circle (1.5pt);

\begin{scope}[shift = {(-3.5,-1.8)}]
\draw[thick] (0,0) -- (1,0);
\draw[thick] (1,0) -- (0.5,0.8660254);
\draw[thick] (0.5,0.8660254) -- (0,0);
\draw[fill = black] (0,0) circle (1.5pt);
\draw[fill = black] (1,0) circle (1.5pt);
\draw[fill = black] (0.5,0.8660254) circle (1.5pt);
\node at (.5,.62) {\tiny{$1$}};
\node at (.5,-.2) {\tiny{$1$}};
\node at (.25,.15) {\tiny{$2$}};
\node at (.87,.5) {\tiny{$2$}};
\node at (.75,.15) {\tiny{$3$}};
\node at (.13,.5) {\tiny{$3$}};
\end{scope}

\begin{scope}[shift = {(-1.5,-1.8)}]
\draw[thick] (0,0) -- (1,0);
\draw[thick] (1,0) -- (0.5,0.8660254);
\draw[thick] (0.5,0.8660254) -- (0,0);
\draw[fill = black] (0,0) circle (1.5pt);
\draw[fill = black] (1,0) circle (1.5pt);
\draw[fill = black] (0.5,0.8660254) circle (1.5pt);
\node at (.5,.62) {\tiny{$1$}};
\node at (.5,-.2) {\tiny{$1$}};
\node at (.25,.15) {\tiny{$2$}};
\node at (.87,.5) {\tiny{$2$}};
\node at (.75,.15) {\tiny{$3$}};
\node at (.13,.5) {\tiny{$3$}};
\end{scope}

\begin{scope}[shift = {(.5,-1.8)}]
\draw[thick] (0,0) -- (1,0);
\draw[thick] (1,0) -- (0.5,0.8660254);
\draw[thick] (0.5,0.8660254) -- (0,0);
\draw[fill = black] (0,0) circle (1.5pt);
\draw[fill = black] (1,0) circle (1.5pt);
\draw[fill = black] (0.5,0.8660254) circle (1.5pt);
\node at (.5,.62) {\tiny{$1$}};
\node at (.5,-.2) {\tiny{$1$}};
\node at (.25,.15) {\tiny{$2$}};
\node at (.87,.5) {\tiny{$2$}};
\node at (.75,.15) {\tiny{$3$}};
\node at (.13,.5) {\tiny{$3$}};
\end{scope}

\begin{scope}[shift = {(2.5,-1.8)}]
\draw[thick] (0,0) -- (1,0);
\draw[thick] (1,0) -- (0.5,0.8660254);
\draw[thick] (0.5,0.8660254) -- (0,0);
\draw[fill = black] (0,0) circle (1.5pt);
\draw[fill = black] (1,0) circle (1.5pt);
\draw[fill = black] (0.5,0.8660254) circle (1.5pt);
\node at (.5,.62) {\tiny{$1$}};
\node at (.5,-.2) {\tiny{$1$}};
\node at (.25,.15) {\tiny{$2$}};
\node at (.87,.5) {\tiny{$2$}};
\node at (.75,.15) {\tiny{$3$}};
\node at (.13,.5) {\tiny{$3$}};
\end{scope}

\end{scope}

\end{tikzpicture}
    \caption{\cref{ex: distance-4} showing a set of distance-3 cliques in a graph $G$, precolored with color palette $[\Delta+3]$, that does not extend to a total-$(\Delta+3)$-coloring of $G$.}
    \label{fig:distance-4-cliques}
\end{figure}

\begin{example}\label{ex: distance-4}
Let $G$ be the graph obtained by taking $K_{1,4}$ and joining to each leaf, a $K_3$ (by join we mean add all three edges between the leaf and the $K_3$). First note that $\Delta(G)=4$. We precolor the four $K_3$'s identically with colors $[3]$ where $\varphi(v)=\varphi(e)$ iff $e$ is not incident to $v$ (see \cref{fig:distance-4-cliques}). We now show that the precoloring $\varphi$ does not extend to a total-$7$-coloring of $G$. Suppose for contradiction that such an extension exists. Then observe, that any leaf of the $K_{1,4}$ along with its incident edges joining the $K_3$, must receive different colors than $[3]$, say $\{4,5,6,7\}$, since each leaf and incident edge joining the $K_3$ sees every color in $[3]$. But this means the edges of the $K_{1,4}$ cannot be colored with $\{4,5,6,7\}$ which is a contradiction since we cannot only use the colors $[3]$ to edge color $K_{1,4}$.
\end{example}

\section{Proof of \cref{thm: D+t - all cliques}}\label{sec: thm1.8}

We first introduce some notation and terminology that will be helpful for the rest of the paper. Given a graph $G$, we let $V_k$ be the vertices of $G$ with degree $k$, and we call any $v\in V_k$ a \emph{$k$-vertex}. We let $V_{[a,b]} = \bigcup_{i=a}^bV_i$. Similarly, given a planar embedding of $G$, we let $F_k$ be the faces of $G$ with length $k$ and we call such $f\in F_k$ a \emph{$k$-face}.

We prove the following slight strengthening of Theorem \ref{thm: D+t - all cliques}, as discussed in the introduction.

\begin{theorem}\label{thm: 111t}
Let $1\leq t\leq 4$, let $G$ be a planar graph with maximum degree $\Delta\geq 28-t$, and let $H$ be a subgraph of $G$. If $H$ is a set of distance-3 cliques, then any $(\Delta+t)$-total-coloring of $H$ can be extended to $G$.
\end{theorem}

\begin{proof}
Suppose that $G$ is a counterexample for which $t$ is maximum and, subject to that, $|V(G)|+|E(G)|$ is minimum. Then there exists a set $H$ of distance-3 cliques in $G$ such that some total-$(\Delta+t)$-coloring $\varphi_0$ of $H$ cannot be extended to $G$.

We proceed with a series of claims about $G$ and $H$.

\begin{claim}\label{clm: Gprime} If $e\in E(G)\setminus E(H)$, then $\varphi_0$ can be extended to $G-e$.
\end{claim}

\begin{proofc} Let $G' = G-e$. If $\Delta(G')\geq 28-t$, then by minimality of $G$, $\varphi_0$ extends to $G'$. If $\Delta(G')<28-t$, then $\Delta(G) = 28-t$, $\Delta(G') = 28-(t+1)$, and we may consider $\varphi_0$ as a total-$(\Delta(G')+(t+1))$-coloring of $H$. 
If $t\leq 3$ then by maximality of $t$,  $\varphi_0$ extends to total-$(\Delta(G')+(t+1))$-coloring of $G'$, i.e., a total-$(\Delta+t)$-coloring of $G'$. If $t=4$ then we instead apply \cref{PlanarDelta+3} to get the same conclusion.
\end{proofc}

Given $v\in V(G)$, we say $v$ is \textit{high} if $\deg(v)\geq \lceil \frac{\Delta+t}{2}\rceil$ and otherwise $v$ is \textit{low}.

\begin{claim}\label{obs: ALL CLIQUES - deg sum}
Let $uv\in E(G)\setminus E(H)$ such that $u\notin V(H)$ and $\deg(u)\leq \frac{\Delta+t}{2}$. Then $\deg(u)+\deg(v)\geq \Delta+t$. Moreover, 
every edge in $E(G)\setminus E(H)$ has at least one high endpoint.
\end{claim}

\begin{proofc} Note that the second sentence follows from the first, since an edge in $E(G)\setminus E(H)$ cannot have both endpoints in $H$ without violating the distance condition on $H$. To prove the first sentence, we take $u, v$ as given and consider $G' = G-uv$. By Claim \ref{clm: Gprime} we know that $\varphi_0$ can be extended to total-$(\Delta+t)$-coloring of $G'$, say $\varphi$. We will try to extend $\varphi$ to $G$. If $\varphi(u)=\varphi(v)$, then our first step in this extension is to modify $\varphi$ by recoloring $u$ (which can be considered since $u\not\in V(H)$). The recoloring is possible since after removing the color from $u$, the number of colors seen by $u$ in $\varphi$ is at most $2(\deg_G(u)-1)+1$, which is at most $\Delta+t-1$ because $\deg(u)\leq \frac{\Delta+t}{2}$. Once we have $\varphi(u)\neq \varphi(v)$, notice that $uv$ sees at most $\deg_G(u)-1+\deg_G(v)-1+2 = \deg_G(u)+\deg_G(v)$ colors. If $\deg(u)+\deg(v)<\Delta+t$, then there exists a color for $uv$ and $\varphi_0$ can be extended to $G$, contradiction. Hence $\deg(u)+\deg(v)\geq \Delta+1$, as desired.
\end{proofc}

\begin{claim}\label{clm: Vt} The following hold.
\begin{enumerate}
    \item[(1)] For $t\geq 2$, $V_{[1,t-1]}\subseteq V(H)$.
    \item[(2)] $V_t\subseteq V(H)$.
    \item[(3)] If $v\in V_{t+1}\setminus V(H)$, then $N(v)\subseteq V_{\Delta}$.
\end{enumerate}
\end{claim}

\begin{proofc} (1) If $t\geq 2$ and $u\in V_{[1,t-1]}\setminus V(H)$, then any edge $uv$ is in $E(G)\setminus E(H)$. So since $t-1\leq \frac{\Delta+t}{2}$, Claim \ref{obs: ALL CLIQUES - deg sum} tells us that $d(v)>\Delta$, contradiction.

(2) Suppose for contradiction there exists $v\in V_t\setminus V(H)$. Let $N(v) = \{u_1,...,u_t\}$. By \cref{obs: ALL CLIQUES - deg sum}, $u_i \in V_{\Delta}$ for all $i\in [t]$. Let $G' = G-v$.
By Claim \ref{clm: Gprime} we know that $\varphi_0$ can be extended to total $(\Delta+t)$-coloring of $G'$, say $\varphi$. We try to extend $\varphi$ to $G$ by first coloring $vu_1, vu_2,..., vu_t$ and finally coloring $v$. Note that by coloring in this sequence when we come to color edge $vu_i$ the number of forbidden colors is at most $\Delta+i-1 < \Delta+t$ colors. When we come to $v$ the number of forbidden colors is at most $2\deg(v) = 2t$. Since $t\leq 4$ and $\Delta\geq 24$, we get that $2t< \Delta+t$ and hence that $\varphi$ can be extended to $G$, contradiction.

(3) Let $v\in V_{t+1}\setminus V(H)$ and $N(v) = \{u_1,...,u_{t+1}\}$.
By \cref{obs: ALL CLIQUES - deg sum}  we get that $N(v)\subseteq V_{[\Delta-1,\Delta]}$. Suppose for contradiction  $u_{t+1}\in V_{\Delta-1}$. By the same argument as in (2), we extend $\varphi_0$ to $G-v$; we then attempt to extend to $G$ by first coloring $vu_1, vu_2,..., vu_t$ and finally coloring $v$. As in (2), we can color $vu_1, vu_2,..., vu_t$ greedily. But now we also have at most 
$\Delta-1+t$ colors forbidden when we come to color $vu_{t+1}$. 
When we come to $v$ the number of forbidden colors is at most $2\deg_G(v)=2t+2\leq 10$. Since $\Delta\geq 24$, we can indeed extend $\varphi_0$ to $G$, contradiction.
\end{proofc}


\begin{claim}\label{obs: ALL CLIQUES - POT}
   When $t\in \{1,2\}$, $|V_{t+1}\setminus V(H)| < |V_{\Delta}|$.
\end{claim}

\begin{proofc}
Let $t\in \{1,2\}$ and let $B$ be the bipartite subgraph $G$ induced by parts $V_{t+1}\setminus V(H)$ and $V_{\Delta}$. Suppose first that $B$ is acyclic. Then $|E(F)| < |V(F)|$. By Claim \ref{clm: Vt}(3), we know that $|E(B)|\geq (t+1)|V_{t+1}\setminus V(H)|\geq 2|V_{t+1}\setminus V(H)|$. Since $|V(B)|=|V_{t+1}\setminus V(H)|+|V_{\Delta}|$, this implies that  $|V_{t+1}\setminus V(H)| < |V_{\Delta}|$, as desired.

Suppose now, for a contradiction, that $B$ contains some cycle $C$.
Let $G'=G-E(C)$. We know that $\varphi_0$ can be extended to total $(\Delta+t)$-coloring of $G'$, say $\varphi$, by Claim \ref{clm: Gprime}  (extend to all but one edge of $C$, then remove the colors from the others).
We attempt to extend $\varphi$ to $G$ as follows. We begin by uncoloring each vertex in $V(C)\cap V_{t+1}$ (none of which are in $H$, by definition of $B$). The number of colors forbidden for coloring any $uv\in E(G)$ is at most $\deg(u)-2+\deg(v)-2+1\leq \Delta+t-2$ colors. Thus, there exists two choices for each edge. Since even cycles are $2$-edge-choosable, we may color $E(C)$. Now when we go to color some $v\in V_{t+1}\cap V(C)$ the number of forbidden colors is at most $2(t+1)\leq 6$. Since $\Delta\geq 5$, $\Delta+t\geq 7$ and thus, we can color $v$. Hence we have extended $\varphi_0$ to $G$, contradiction. 
\end{proofc}

We now fix a planar embedding of $G$ for the remainder of the proof. The graph $H$ is made up of vertex-disjoint cliques, and we label these cliques as $H_1, \ldots, H_m$. For each $H_i$, $1\leq i\leq m$, we define the \textit{configuration of $H_i$}, denoted $\mathcal{C}(H_i)$, as $G[N[H_i]]$ where for $X\subseteq V(G)$, $N[X] = X \cup N(X)$. Since $H$ is a set of distance-3 cliques, every vertex in $G$ is contained in at most one configuration. If $H_i$ is a $k$-clique with vertices $v_1,...,v_k$ and $\deg(v_1)\leq ...\leq \deg(v_k)$, we say $\mathcal{C}(H_i)$ is a $(\deg(v_1),...,\deg(v_k))$-configuration. We say that a configuration $\mathcal{C}(H_i)$ is \emph{poor} if it contains at most one high vertex.  Otherwise we say it is \emph{rich}. As it turns out, there are only $16$ poor configurations, each depicted in \cref{fig:poor-configurations-1}.

\begin{claim}\label{lem:ALL CLIQUES-enum-poor-conf}
The only poor configurations are the ones depicted in \cref{fig:poor-configurations-1}. Moreover, if $\mathcal{C}(H_i)$ is a $(2,2)$-, $(2,3,3)$-, $(3,3,3)$-, $(3,3,4,4)$-, or $(3,4,4,4)$-configuration that is poor, then all the vertices in $V(H_i)$ of degree $|V(H_i)|$ are adjacent to a unique high vertex in $\mathcal{C}(H_i)$.
\end{claim}

\begin{proofc}
We see from Figure \ref{fig:poor-configurations-1} that the second sentence follows from the first; we must prove the first sentence of the claim. Let $\mathcal{C}(H_i)$ be a poor configuration, with $h$ the degree of the unique high vertex, if it exists. We start by showing that all of the following properties hold.
\begin{enumerate}
    \item[(i)] The number of high vertices in $V(H_i)$ is at most $1$.
    \item[(ii)] If the number of high vertices in $V(H_i)$ is $1$, then all other vertices in $V(H_i)$ must have degree $|V(H_i)|-1$.
    \item[(iii)] The degree of any low vertex in $V(H_i)$ is at most $|V(H_i)|$
    \item[(iv)] All $|V(H_i)|$-vertices in $V(H_i)$ share a common high neighbor.
\end{enumerate} 

Property (i) follows from the definition. For (ii), if $v\in V(H_i)$ $v$ low and $\deg(v)\geq |V(H_i)|$, then there must exists some $v'\in N(v)\setminus V(H_i)$ which implies that $vv'\notin E(H)$. Thus, by \cref{obs: ALL CLIQUES - deg sum}, $v'$ is high and hence $\mathcal{C}(H_i)$ is rich, contradiction. 
Similarly, for (iii), if $V(H_i)$ contains a low vertex $v$ with $\deg(v)\geq |V(H_i)|+1$, then there must exist two distinct neighbors of $v$, say $v'$ and $v''$ with $vv',vv''\notin E(H)$. Thus, by \cref{obs: ALL CLIQUES - deg sum} $v'$ and $v''$ are high and hence $\mathcal{C}(H_i)$ is rich, contradiction. By the same argument again, if $V(H_i)$ contains vertices which have degree $|V(H_i)|$, then they must be adjacent to a shared high vertex as desired in (iv).

First suppose $|V(H_i)|=1$ with $V(H_i)=\{v\}$.  Property (iii) tells us that $v$ must be high (since $G\not \cong K_1$). Thus, $\mathcal{C}(H_i)$ is an $(h)$-configuration, as depicted in the top-left of Figure \ref{fig:poor-configurations-1}. 

Next suppose $|V(H_i)|=2$ with $V(H_i)=\{u,v\}$ and $\deg(u)\leq \deg(v)$. By (i), $u$ must be low and then by (iii), $\deg(u)\leq 2$. If $\deg(u)=1$, then $\deg(v)\geq 2$ since $G\not \cong K_2$. By (iii), this means that when $\deg(u)=1$, the vertex $v$ is either a $2$-vertex or a high vertex, and either way we have a $(1,2)$- or $(1,h)$-configuration, as depicted in the top center of Figure \ref{fig:poor-configurations-1}. If $\deg(u)=2$, then by (ii) and (iii), $\deg(v)=2$ and we have a $(2,2)$-configuration, as pictured in the top-right of Figure \ref{fig:poor-configurations-1}.

Now suppose $|V(H_i)|=3$ with $V(H_i)=\{u,v,w\}$ and $2\leq \deg(u)\leq \deg(v)\leq \deg(w)$. By (i), $u$ and $v$ must be low and hence by (iii) $\deg(u)\leq \deg(v)\leq  3$. If $w$ is high then by (ii), $\deg(u) = \deg(v) = 2$ and if $w$ is low then by (iii), $\deg(w) = 3$ (since $G\not \cong K_3$,). Thus we have either a $(2,2,h)$-, $(2,2,3)$-, $(2,3,3)$- or a $(3,3,3)$-configuration, as pictured in the second row of Figure \ref{fig:poor-configurations-1}. 

Finally, suppose $|V(H_i)|=4$ with $V(H_i)=\{u,v,w,x\}$ and $3\leq \deg(u)\leq \deg(v)\leq \deg(w)\leq \deg(x)$. By (i), $u,v,w$ are low and hence by (iii), $3\leq \deg(u)\leq \deg(v)\leq \deg(w)\leq 4$. If $x$ is high, then by (ii), $\deg(u)=\deg(v)=\deg(w)=3$ and we have a $(3,3,3,h)$-configuration as depicted in the bottom row of Figure \ref{fig:poor-configurations-1}. If $x$ is low, then by (iii), and since $G\not \cong K_4$, $\deg(x)=4$. We cannot have a $(4,4,4,4)$-configuration since by (iv), all $4$-vertices would share a common neighbor which would mean $G\cong K_5$, contradiction. Thus we have either a $(3,3,3,4)$- $(3,3,4,4)$-, $(3,4,4,4)$-configuration, as depicted in the bottom row of Figure \ref{fig:poor-configurations-1}. 
\end{proofc}

\begin{figure}[h!]
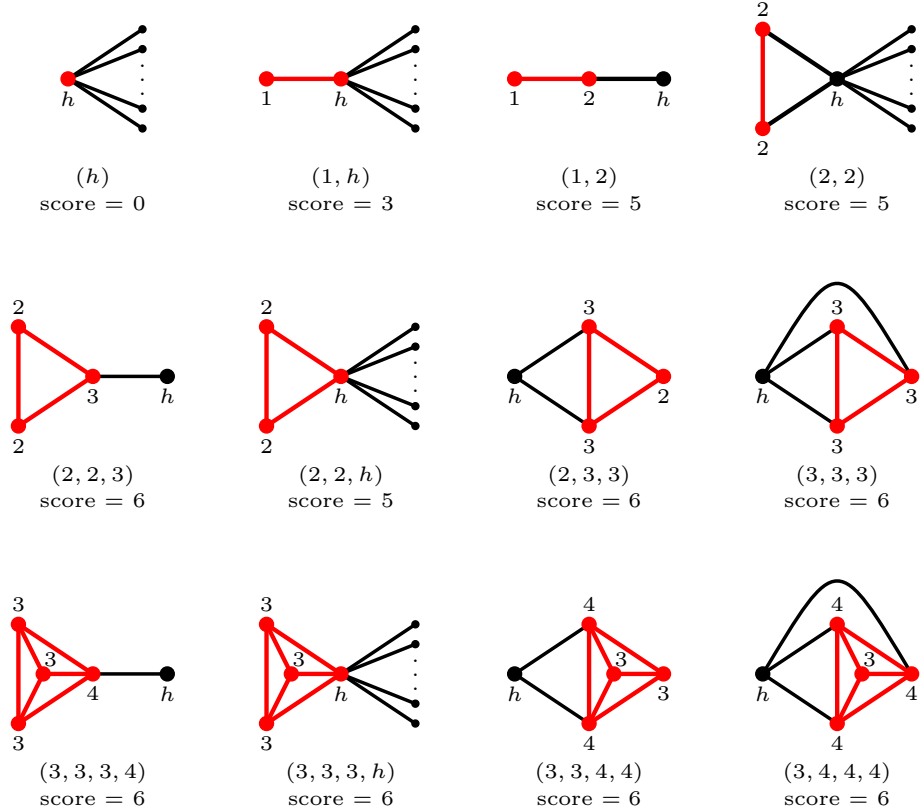

    \centering
    \includestandalone[width=0.75\textwidth]{poor-configurations-all-cliques}
    \caption{All possible poor configurations are depicted (with black edges) along with their corresponding cliques (with red edges and vertices). Vertex degrees are indicated for all the clique vertices, and $h$ indicates the degree of the high vertex, which exists in each case.}
    \label{fig:poor-configurations-1}
\end{figure}

We define a \emph{score}  for each $H_i$, $1\leq i\leq m$, as follows:
\[s(H_i) = \sum_{v\in V_{[1,3]}\cap V(H_i)}4-\deg(v)+|\{f\in F_3: \text{$f$ contains at least one precolored edge from $H_i$}\}|.\]
Note that the score of each $H_i$ is related to the amount of charge that it will somehow need in our discharging argument, which is still to come. The following claim establishes an upper bound for the score of every $H_i$.

\begin{claim}\label{lem: ALL CLIQUES - SCORE}
For any $H_i$, $s(H_i)\leq 6$ and moreover if $\mathcal{C}(H_i)$ is poor, then \[s(H_i)\leq \begin{cases}
0 \hspace{.5cm} &\text{if $|V(H_i)|=1$}\\
5  &\text{if $|V(H_i)|=2$}\\
6  &\text{if $|V(H_i)|\geq3$}\\
\end{cases}\]
\end{claim}

\begin{proofc}
By \cref{lem:ALL CLIQUES-enum-poor-conf}, the moreover statement can be easily verified by checking the configurations in \cref{fig:poor-configurations-1}.
We now show that for any configuration $\mathcal{C}(H_i)$, $s(H_i)\leq 6$.

First suppose $|V(H_i)|=1$ with $V(H_i)=\{v\}$. Then the only thing contributing to $s(H_i)$ is the $4-\deg(v)$ if $v\in V_{[1,3]}$. Thus, $s(H_i)\leq 3$.

Next suppose $|V(H_i)|=2$ with  $V(H_i)=\{u,v\}$. If either $u$ or $v$ is a 1-vertex (note both cannot be) then there is no $3$-face incident to the precolored edge $uv$ and hence $s(H_i)\leq 3+2=5$. If neither $u$ nor $v$ is a $1$-vertex, then there are at most two $3$-faces incident to the precolored edge $uv$ and thus, $s(H_i)\leq 2+2+2=6$.

Next suppose $|V(H_i)|=3$ with $V(H_i)=\{u,v,w\}$. If two of these vertices are $2$-vertices, then the only possible $3$-face incident with $uv,uw$, or $vw$ is the $3$-clique itself, and thus, $s(H_i)\leq 2+2+1+1=6$. If exactly one vertex in $\{u,v,w\}$ is a $2$-vertex, say without loss of generality $u$, then there are at most two $3$-faces incident to an edge in $uv, uw, vw$, namely the $3$-clique itself and possibly a $3$-face incident with $vw$. Thus, $s(H_i)\leq 2+1+1+2=6$. We may now assume no vertex in $\{u,v,w\}$ is a $2$-vertex. If all three are $3$-vertices, then since $G\not\cong K_4$, there are at most three $3$-faces incident to the edges $uv, uw, vw$ and thus, $s(H_i)\leq 1+1+1+3=6$. Finally suppose there is at least one vertex in $\{u,v,w\}$ with degree at least $4$, then there are at most four $3$-faces incident to an edge in $uv, uw, vw$, and $s(H_i)\leq 1+1+4=6$.

Finally suppose $|V(H_i)|=4$, with vertices $u,v,w,x$. Since $G\not\cong K_4$, at least one vertex must have degree at least $4$. If $k$ vertices in $\{u,v,w,x\}$ are $3$-vertices for some $k\in \{0,1,2,3\}$, then there are at most $6-k$ $3$-faces incident to an edge in $uv, uw, ux, vw, vx, wx$ and thus, $s(H_i)\leq k+(6-k)=6$. \end{proofc}

Consider now a poor configuration $\mathcal{C}(H_i)$ with $|V(H_i)|\geq 2$. We know that the configuration has a high vertex, say $v$, that it is a cut-vertex in $G$, and moreover, that $H_i-v$ is a component of $G-v$ (see Figure \ref{fig:poor-configurations-1}). Thus, $V(H_i)\setminus \{v\}$ lies in some unique face $f'$ of $G-(V(H_i)\setminus \{v\})$. We say $f\in F(G)$ is a \emph{helpful face with respect to $H_i$} if $f$ can be obtained from by extending $f'$ to include at least one edge of $H_i$. It is clear that every $H_i$ has a unique helpful face, but a face can be helpful with respect to multiple cliques $H_i$ (see Figure \ref{fig:helpfulface}). Even so, in the following claim we show that the length $\ell(f)$ of such a face $f$ is so long that it is still going to be useful in our upcoming discharging argument.

\begin{figure}[h!]
    \centering
    \begin{tikzpicture}

\def\R{70pt}
\draw[fill=white] (0,0) circle (\R);

\fill[black] ({\R*cos(70)},{\R*sin(70)}) circle (2pt);
\fill[black] ({\R*cos(110)},{\R*sin(110)}) circle (2pt);
\fill[black] ({\R*cos(-70)},{\R*sin(-70)}) circle (2pt);
\fill[black] ({\R*cos(-110)},{\R*sin(-110)}) circle (2pt);

\begin{scope}[shift={(90:\R)}, rotate=90]
  \draw[line width=1.25pt,red] (-.75,0) -- (0,0);
  \draw[red,fill=red] (-.75,0) circle (2pt);
  \draw[red,fill=red] (0,0) circle (2pt);
  \node at (0.25,0) {\tiny{$h$}};
\end{scope}

\begin{scope}[shift={(150:\R)}, rotate=150]
  \draw[line width=1.25pt,red] (-.75,-.5) -- (-.75,.5);
  \draw[line width=1.25pt,black] (-.75,-.5) -- (0,0);
  \draw[line width=1.25pt,black] (0,0) -- (-.75,.5);
  \draw[red,fill=red] (-.75,-.5) circle (2pt);
  \draw[red,fill=red] (-.75,.5) circle (2pt);
  \draw[black,fill=black] (0,0) circle (2pt);
  \node at (0.25,-0.1) {\tiny{$h$}};
\end{scope}

\begin{scope}[shift={(210:\R)}, rotate=210]
  \draw[line width=1.25pt,red] (-1.5,-.5) -- (-1.5,.5);
  \draw[line width=1.25pt,red] (-1.5,-.5) -- (-0.75,0);
  \draw[line width=1.25pt,red] (-0.75,0) -- (-1.5,.5);
  \draw[line width=1pt] (0,0) -- (-0.75,0);
  \draw[red,fill=red] (-1.5,-.5) circle (2pt);
  \draw[red,fill=red] (-1.5,.5) circle (2pt);
  \draw[red,fill=red] (-0.75,0) circle (2pt);
  \draw[fill=black] (0,0) circle (2pt);
  \node at (0.25,0) {\tiny{$h$}};
\end{scope}

\begin{scope}[shift={(270:\R)}, rotate=270]
  \draw[line width=1.25pt,red] (-.75,-.5) -- (-.75,.5);
  \draw[line width=1.25pt,red] (-.75,-.5) -- (0,0);
  \draw[line width=1.25pt,red] (0,0) -- (-.75,.5);
  \draw[red,fill=red] (-.75,-.5) circle (2pt);
  \draw[red,fill=red] (-.75,.5) circle (2pt);
  \draw[red,fill=red] (0,0) circle (2pt);
  \node at (0.2,0) {\tiny{$h$}};
\end{scope}

\begin{scope}[shift={(330:\R)}, rotate=150]
  \draw[line width=1.25pt,red] (0.75,-.5) -- (0.75,.5);
  \draw[line width=1.25pt,red] (0.75,-.5) -- (1.5,0);
  \draw[line width=1.25pt,red] (1.5,0) -- (0.75,.5);
  \draw[line width=1pt, bend left=0] (0,0) to (0.75,.5);
  \draw[line width=1pt, bend right=0] (0,0) to (0.75,-.5);
  \draw[line width=1pt] (0,0) .. controls (0.75,1.25) .. (1.5,0);
  \draw[red,fill=red] (0.75,-.5) circle (2pt);
  \draw[red,fill=red] (0.75,.5) circle (2pt);
  \draw[red,fill=red] (1.5,0) circle (2pt);
  \draw[fill=black] (0,0) circle (2pt);
  \node at (-0.25,0) {\tiny{$h$}};
\end{scope}

\begin{scope}[shift={(30:\R)}, rotate=210]
  \draw[line width=1.25pt,red] (0.75,-.5) -- (0.75,.5);
  \draw[line width=1.25pt,red] (0.75,-.5) -- (1.5,0);
  \draw[line width=1.25pt,red] (1.5,0) -- (0.75,.5);
  \draw[line width=1.25pt,red] (1.0,0) -- (0.75,.5);
  \draw[line width=1.25pt,red] (0.75,-.5) -- (1.0,0);
  \draw[line width=1.25pt,red] (1.5,0) -- (1.0,0);
  \draw[line width=1pt, bend left=0] (0,0) to (0.75,.5);
  \draw[line width=1pt, bend right=0] (0,0) to (0.75,-.5);
  \draw[red,fill=red] (0.75,-.5) circle (2pt);
  \draw[red,fill=red] (0.75,.5) circle (2pt);
  \draw[red,fill=red] (1.5,0) circle (2pt);
  \draw[red,fill=red] (1.0,0) circle (2pt);
  \draw[fill=black] (0,0) circle (2pt);
  \node at (-0.25,0) {\tiny{$h$}};
\end{scope}

\node at (0,0) {\large{$f$}};

\end{tikzpicture}
    \caption{An example of a face $f$ which is the helpful face to $6$ different poor configurations. }
    \label{fig:helpfulface}
\end{figure}
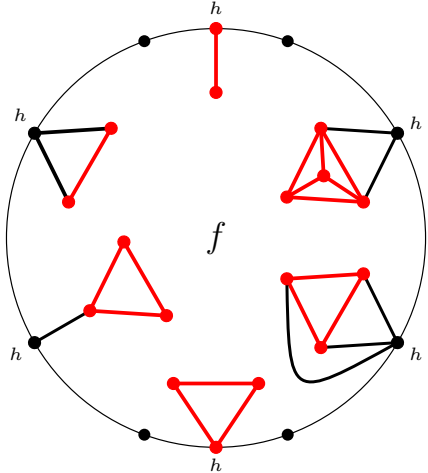

\begin{claim}\label{lem: ALL CLIQUES - FACE LENGTH}
 Let $f\in F(G)$ and suppose $f$ is the helpful face of $x_2+x_3$ different poor configurations, with $x_2$ being the number of these poor-configurations $\mathcal{C}(H_i)$ with $|V(H_i)|=2$, and $x_3$ being the number of these poor configurations $\mathcal{C}(H_j)$ with $|V(H_j)|\geq 3$. Then $$\ell(f)\geq 4(x_2+x_3)$$ and moreover if $x_2+x_3=1$, then \[\ell(f)\geq \begin{cases}
     5 \hspace{.5cm} &\text{if $x_2=1$ and $x_3=0$}\\
     6 &\text{if $x_2=0$ and $x_3=1$}.\\
 \end{cases}\]
\end{claim}

\begin{proofc}
Let $W$ be the boundary walk of $f$. If $W$ does not contain cycle then clearly $\ell(f)\geq 4(x_1+x_2)$ since each poor configuration has at least $2$ edges contained in $W$ which are each counted twice in a face without a cycle. Now suppose $W$ contains an cycle $C = v_1,...,v_k$. Note that $k\geq \max\{x_2+x_3,3\}$. 

For each configuration counted by $x_3$, note there is an additional length of at least $3$ in $W$ which is not accounted for in $C$. So overall, we know that $\ell(f)\geq k+3x_3$.

For each configuration counted in $x_2$, there is an additional length of at least $2$ in $W$ that not counted by $C$. However, looking at Figure \ref{fig:poor-configurations-1} we see that equality is only achieved when the poor configuration is a $(1,h)$-configuration and the high vertex (which is precolored) lies on $C$ (see the top configuration of \cref{fig:helpfulface}). In this case however, if $v_i$ is the high vertex of the configuation, then since precolored cliques are at a distance of $3$, we get that $v_{i-1}$ and $v_{i+1}$ are not precolored and hence not the high vertex of a poor configuration. So in fact $C$ itself must have length at least one longer to compensate for each such configuration, and the configuration still adds at least 3 to $\ell(f)$. 

We have now seen that in general, we have $\ell(f)\geq (x_1+x_2)+ 3x_1+3x_2 = 4(x_1+x_2).$ If $x_1+x_2=1$, then $\ell(f)\geq 3+2x_2= 5$ when $x_2=1,x_3=0$ and $\ell(f)\geq 3+3x_2 = 6$ when $x_2=0,x_3=1$. 
\end{proofc}

We are now finally ready to define an initial charge $\mu$ on $x\in V(G)\cup F(G)$, on every configuration $\mathcal{C}(H_i)$ in $G$, and on a \emph{global pot} $\mathcal{P}$. Initially, let $\mu(\mathcal{C}(H_i))=\mu(\mathcal{P})=0$ for every configuration in $G$ and let \[\mu(x) = \begin{cases}
    \deg(x)-4 \hspace{.5cm}&\text{if $x\in V(G)$}\\
    \ell(x)-4 &\text{if $x\in F(G)$}\\
\end{cases}\]
By Euler's formula we see that \[\mu(G) = \sum_{H_i\in H}\mu(\mathcal{C}(H_i))+\sum_{v\in V(G)}\mu(v)+\sum_{f\in F(G)}\mu(f) +\mu(\mathcal{P}) = 4|E(G)|-4|V(G)|-4|F(G)| = -8.\]
The items in $V(G)\cup F(G)$ which receive a negative initial charge are $3$-faces, and vertices in $V_{[1,3]}$. We now present the discharging rules, which we perform in order.

\bigskip    

\begin{mdframed}
    
\textit{Discharging rules}

\begin{enumerate}[label=(R\arabic*)]
\item \label{R1} Each rich configuration $\mathcal{C}(H_i)$ receives $\frac{s(H_i)}2$ charge from each high vertex in $\mathcal{C}(T)$.
\item \label{R2} Each poor configuration $\mathcal{C}(H_i)$ with $|V(H_i)|=2$ receives charge $1$ from its helpful face and receives $s(H_i)-1$ from its high vertex. Each poor configuration $\mathcal{C}(H_i)$ with $|V(H_i)|\geq 3$ receives charge $2$ from its helpful face and receives $s(H_i)-2$ from its high vertex.
\item \label{R3} Each $v\in V_{[1,3]}\cap V(H)$ receives charge $4-\deg(v)$ from its corresponding configuration.
\item \label{R4} Each face $f\in F_3$ with at least one edge in $E(H)$ receives charge $1$ from its corresponding configuration.
\item \label{R5} Each face $f\in F_3$ with no edges in $E(H)$ receives charge $\frac{1}{2}$ from its incident high vertices.
\item\label{R6} When $t\leq 2$, each vertex $u\in V_{\Delta}$ gives charge $4-\deg(v)$ to $\mathcal{P}$ and each vertex $v\in V_{t+1}\setminus V(H)$ receives charge $4-\deg(v)$ from $\mathcal{P}$.
\item\label{R7} When $t=1$, Each vertex $v\in V_3\setminus V(H)$ receives charge $\frac{1}{3}$ from each neighbor.
\end{enumerate}
\end{mdframed}

\medskip

Let $\mu'$ be the final charge after \ref{R1}-\ref{R7} have been applied. We'll prove $\mu'\geq 0$, contradicting  $\mu'(G) = \mu(G) = -8$ and completing the proof. The rest of the proof consists of one claim each for $\mathcal{P}$, faces, configurations and vertices, respectively.

\begin{claim}
$\mathcal{P}$ has nonnegative final charge.
\end{claim}

\begin{proofc}
 The only rule affecting $\mathcal{P}$ is \ref{R6}. By \cref{obs: ALL CLIQUES - POT}, $\mu'(\mathcal{P}) \geq \mu(\mathcal{P})+|V_{[\Delta-1,\Delta]}|-|V_{t+1}\setminus V(H)| > 0$.  
\end{proofc}

\begin{claim}
Every face has nonnegative final charge.
\end{claim}

\begin{proofc}
Let $f\in F(G)$. If $f$ is neither a $3$-face, nor is the helpful face for any configuration, then $f$ is not affected by the discharging rules and $\mu'(f)=\mu(f)=\ell(f)-4\geq 4-4=0$.

Now suppose $f\in F_3$. Since any helpful face must have length at least $5$ by \cref{lem: ALL CLIQUES - FACE LENGTH}, $f$ is only affected by one of \ref{R4}, \ref{R5}. In the former case, $\mu'(f)=\mu(f)+1=3-4+1=0$. In the latter case, then by \cref{obs: ALL CLIQUES - deg sum}, $f$ must be incident to at least $2$ high vertices. Thus, by \ref{R5} $\mu'(f)=\mu
(f)+2\cdot \frac{1}{2} = 3-4+1=0$.

Finally suppose $f$ is the helpful face for $x_2+x_3$ poor configurations, with $x_2, x_3$ as in the statement of Claim \ref{lem: ALL CLIQUES - FACE LENGTH}. By \ref{R2}, $\mu'(f) = \mu(f) - x_1-2x_2 = \ell(f)-4-x_1-2x_2$. If $x_2+x_3\geq 2$, then by \cref{lem: ALL CLIQUES - FACE LENGTH}, $\mu'(f)\geq 4x_2+4x_3-4-x_2-2x_3=3x_2+2x_3-4\geq 0$. If $x_2+x_3=1$ then by \cref{lem: ALL CLIQUES - FACE LENGTH} $\mu'(f)\geq 5-4-1=0$ when $x_2=1,x_3=0$ and $\mu'(f)\geq 6-4-2=0$ when $x_2=0,x_3=1$. In all cases, $\mu'(f)\geq 0$.
\end{proofc}

\begin{claim}
    Each configuration has nonnegative final charge. 
\end{claim}

\begin{proofc}
    Let $\mathcal{C}(H_i)$ be a configuration. By \ref{R1} and \ref{R2} $\mathcal{C}(H_i)$ receives $s(H_i)$ charge. Then by \ref{R3} and \ref{R4} $\mathcal{C}(H_i)$ gives $2|V_2\cap V(H_i)|+|V_3\cap V(H_i)|+|\{f\in F_3: E(f)\cap E(H)\neq \emptyset\}| = s(H_i)$. Thus $\mu'(\mathcal{C}(H_i)) = \mu(\mathcal{C}(H_i))+s(H_i)-s(H_i)= 0$.
\end{proofc}

\begin{claim}
Every vertex has nonnegative final charge.
\end{claim}

\begin{proofc}
Let $v\in V(G)$. First suppose $v\in V_{1}$ or $v\in V_{[2,3]}\cap V(H)$. Then by \cref{clm: Vt}  $v\in V(H)$, and thus, it is contained in a configuration. By \ref{R3}, $\mu'(v) = \mu(v)+ 4-\deg(v)=\deg(v)-4+4-\deg(v)=0$. Next suppose $v\in V_{[2,3]}\setminus V(H)$. Then by \cref{clm: Vt} we know $t\leq 2$. If $v\in V_{t+1}$, then by \ref{R6}, $\mu'(v)=\mu(v)+4-\deg(v)=\deg(v)-4+4-\deg(v)=0$. If $v\in V_{t+2}$, then $t=1$ and $v\in V_3$. By \ref{R7} $\mu'(v)= \mu(v)+\frac{1}{3}\deg(v)=-1+1=0$.
We may now assume $\deg(v)\geq 4$. 

If $v$ is low then it is not affected by any of the discharging rules and thus, $\mu'(v)=\mu(v)=\deg(v)-4\geq 0$. We may therefore assume that $v$ is high. Note that since $\Delta\geq 28-t$, $\deg(v)\geq \lceil \frac{\Delta+t}{2} \rceil\geq 14$ which will turn out to be exactly what we need.

Since there are extra rules for vertices in $V_{\Delta}$ when $t=2$ and for $V_{[\Delta-2,\Delta]}$ when $t= 1$, we take care of those separately later. So suppose $v$ is high but if $t= 2$ then $v\notin V_{\Delta}$ and if $t=1$ $v\notin V_{[\Delta-2,\Delta]}$.

Suppose first that $v$ is in a poor configuration $\mathcal{C}(H_i)$. If $|V(H_i)|=1$, then by \ref{R5}, $\mu'(v)\geq \mu(v)-\frac{1}{2}\deg(v) = \frac{1}{2}\deg(v)-4$. Thus, $\mu'(v)\geq 0$ because $\deg(v)\geq 8$. 
If $|V(H_i)|=2$ then $v$ gives charge $s(H_i)-1$ to its poor configuration by \ref{R2}. The only other rule which may affect $v$ is \ref{R5} where it gives charge $\frac{1}{2}$ to any incident $3$-face with no edges in $E(H)$. Since $v$ is a cut vertex in $G$, the number of incident faces to $v$ is at most $\deg(v)-1$. Moreover, at least one of these faces cannot be a $3$-face with no edges in $E(H)$ (since there is at least one face incident to $v$ and $H_i$). Thus, $\mu'(v)\geq \mu(v)-(s(H_i)-1+\frac{1}{2}(\deg(v)-2))$. By \cref{lem: ALL CLIQUES - SCORE}, $\mu'(v)\geq \deg(v)-4-(5-1+\frac{1}{2}(\deg(v)-2))=\frac{1}{2}\deg(v)-7$. Thus, $\mu'(v)\geq 0$ because $\deg(v)\geq 14$. Similarly, if $|V(H_i)|\geq 3$, then by \ref{R2} and \ref{R5}, $\mu'(v)\geq \mu(v)-(s(H_i)-2+\frac{1}{2}(\deg(v)-2))$. By \cref{lem: ALL CLIQUES - SCORE}, $\mu'(v)\geq \deg(v)-4-(6-2+\frac{1}{2}(\deg(v)-2))=\frac{1}{2}\deg(v)-7$ and hence $\mu'(v)\geq 0$ because $\deg(v)\geq 14$. 

Suppose now that $v$ is in a rich configuration or $v$ is in no configuration at all (and again if $t = 2$, then $v\notin V_{\Delta}$ and if $t=1$ then $v\notin V_{[\Delta-2,\Delta]}$). Then, by \ref{R1} and \ref{R5}, $\mu'(v)\geq \mu(v)-(\frac{s(H_i)}{2}+\frac{1}{2}\deg(v))$. By \cref{lem: ALL CLIQUES - SCORE}, $\mu'(v)\geq \deg(v)-4-(3+\frac{1}{2}\deg(v))=\frac{1}{2}\deg(v)-7$. Since $\deg(v)\geq 14$ we get $\mu'(v)\geq 0$, as desired.

Finally we consider the cases when $t= 2$ and $V_{\Delta}$ or $t=1$ and $v\in V_{[\Delta-2,\Delta]}$. First suppose $t=2$ and $v\in V_{\Delta}$. The only difference between the previous counts is that $v$ also gives charge $1$ by \ref{R6}. Thus, $\mu'(v)\geq \frac{1}{2}\deg(v)-7-1=\frac{1}{2}\Delta-8\geq 0$ since $\Delta\geq 16$. Lastly suppose $t=1$ and $v\in V_{[\Delta-2,\Delta]}$. If $v$ has two $3$-neighbors $u,u'\notin V(H)$, then $uu'\notin E(G)$ and hence there cannot be a $3$-face incident to $v$ and $u$ and $u'$. Thus, if $v$ is adjacent to $s$ vertices in $V_3\setminus V(H)$ and incident to $s'$ $3$-faces which contain no precolored edge, then $s\leq \deg(v), s'\leq \deg(v)$, and $s+s'\leq \frac{3}{2}\deg(v)$. By \ref{R7} $v$ has to give charge $\frac{1}{3}$ to each $3$-neighbor in $V_3\setminus V(H)$ and as mentioned before, $v$ gives charge $\frac{1}{2}$ to each $3$-face with no precolored edge. Also, if $v\in V_{\Delta}$ then by \ref{R6} $v$ gives charge $2$ to $\mathcal{P}$ instead of charge $1$ seen before. Thus, we can bound the previous count of final charge by subtracting an additional $1$ for \ref{R6} if $\deg(v)= \Delta$ and an additional $\frac{1}{3}\cdot ( \frac{1}{2}(\deg(v)))$ if $v\in V_{[\Delta-2,\Delta]}$. If $v\in V_{[\Delta-2,\Delta-1]}$ we see that $\mu'(v)\geq \frac{1}{2}\deg(v)-8-\frac{1}{3}(\frac{1}{2}\deg(v)) = \frac{1}{3}\deg(v)-8 = \frac{1}{3}(\Delta-2)-8$. thus, $\mu'(v)\geq 0$ since $\Delta\geq 26$. If $v\in V_{\Delta}$, then we have $\mu'(v)\geq \frac{1}{3}\deg(v)-9 =\frac{1}{3}\Delta-9\geq 0$ since $\Delta\geq 27$.
\end{proofc}
\end{proof}

\section{Proof of \cref{thm:Planar-Delta+t}}\label{sec: thm1.6}

Before proving \cref{thm:Planar-Delta+t} in general we prove it when $G$ is planar and bipartite; in fact we'll prove a stronger result in this special case. 

\begin{theorem}[Alon and Tarsi]\label{thm:ALON-TARSI} Every planar bipartite graph is 3-choosable.
\end{theorem}

\begin{theorem}[Borodin, Kostochka, Woodall]\label{thm: BORODIC-KOSTOCHKA}  Let $G$ be a bipartite graph and let
  $L$ be an edge list assignment on $G$. If
  $|{L(xy)}| \geq \max\{\deg(x), \deg(y)\}$ for every edge
  $xy \in E(G)$, then $G$ is $L$-edge-colorable.
\end{theorem}

\begin{lemma}\label{lem: planar-bip-1} Let $G$ be a planar bipartite graph and let $L$ be an total list assignment on $G$. If $|L(v)|\geq 3$ for all $v\in V(G)$ and
  $|{L(xy)}| \geq \max\{\deg(x)+2, \deg(y)+2\}$ for every edge
  $xy \in E(G)$, then $G$ is $L$-totally-colorable.
\end{lemma}

\begin{proof} First give an $L$-coloring of the vertices of $G$, which is possible by \cref{thm:ALON-TARSI}. For every edge $e$, modify $L(e)$ to $L'(e)$ by removing the colors of the incident vertices, if they appear. So $|L'(e)| \geq |L(e)|-2$. But then for $e=xy$,  $|{L(xy)}| \geq \max\{\deg(x), \deg(y)\}$. Hence by \cref{thm: BORODIC-KOSTOCHKA}, the edges of $G$ can all be $L'$-colored, giving an $L$-total-coloring of $G$.
\end{proof}


\begin{theorem}\label{thm: planar-bip-prelim} Let $G$ be a planar bipartite graph, and let $H$ be a subgraph of $G$. If $H$ has maximum degree at most $d$, then any $(\Delta+d+4)$-total-coloring of $H$ in $G$ can be extended to $G$.
\end{theorem}

\begin{proof}
Let $G' = G-E(H)$. For each edge $e$ and vertex $v$ in $G'$, let $L'(e),L'(v)$ be obtained from $[\Delta+d+4]$ by removing all colors used on the vertices and edges of $H$ incident and adjacent to $e$ and $v$. Let $xy$ be an arbitrary edge of $G'$. Note that $|L'(xy)|\geq \Delta+d+4-\deg_H(x)-\deg_H(y)-2$. This means that  \[|L'(xy)|\geq \max\{\Delta+2 -\deg_H(x), \Delta+2-\deg_H(y)\}=\max\{\deg_{G'}(x)+2, \deg_{G'}(y)+2\}.\] 

Now let $v$ be an arbitrary vertex of $G'$ that is not colored under $H$. Then \[|L'(v)|\geq \Delta+d+4-2\deg_H(v) \geq 4.\]

Our result now follows by \cref{lem: planar-bip-1}.
\end{proof}


We now prove \cref{thm:Planar-Delta+t} which we restate here for convenience.

\begin{theorem*}[\ref{thm:Planar-Delta+t}] Let $G$ be a planar graph with maximum degree $\Delta$ and let $H$ be a subgraph of $G$. If $H$ has maximum degree $d$, then any $(\Delta+d+6)$-total-coloring of $H$ in $G$ can be extended to $G$. 
\end{theorem*}

\begin{proof}
If $H=\emptyset$ then since planar graphs are $4$-colorable and $(\Delta+1)$-edge-colorable, we trivially get a total-$(\Delta+5)$-coloring of $G$.

For any fixed $\Delta$ and $d$, we choose a counterexample
$(G, H)$ where the quantity \[q(G)=3|{E(G)}| + |V_{[2,d+5]}|\]
is as small as possible. Throughout the proof, we let $\varphi$ be the total-$(\Delta+d+6)$-coloring of $H$. Our goal is to arrive at a contradiction by extending $\varphi$ to $G$.

\begin{claim}\label{2tplus1empty}
  $V_{[2,d+5]}=\emptyset$.
\end{claim}
\begin{proofc}
  Suppose not, and take $v \in V_{[2,d+5]}$. If $v\notin V(H)$, we first color $v$ with any available color and add $v$ to $H$. Since $v$ sees at most $\deg(v)$ colors, such a choice is possible. Note, adding $v$ to $H$ does not increase the maximum degree of $H$ since no edges are added to $H$. We may now assume $v\in V(H)$.

We now form $G'$ from $G$ by deleting $v$, and attaching a leaf $v_u$ to each $u\in N_G(v)$. By construction, $|E(G)|=|E(G')|$ and thus, $q(G')<q(G)$ since the vertex $v$ in $V_{[2,d+5]}$ has been replaced by $\deg(v)$ vertices of degree $1$. For any edge $uv\in E(H)$ we assign $\varphi(uv_u)=\varphi(uv)$ and for all $u\in N_G(v)$, we assign $\varphi(v_u)=\varphi(v)$. Now by minimum counterexample, we can extend $\varphi$ to a total-$(\Delta+d+6)$-coloring of $G'$. Since $\varphi(v_u)=\varphi(v)$ for all $u\in N(v)$, we would like to obtain a total-$(\Delta+d+6)$-coloring of $G$ by merging all such $v_u$ back into one vertex $v$. Note that the resulting coloring may not be proper if $\varphi(uv_u)=\varphi(u'v_{u'})$ for two $u,u'\in N(v)$. If this occurs, we know at most one of $uv$ or $u'v$ can be in $E(H)$ since $H$ is a proper total-coloring. For any such edge $uv$ which is \emph{not} in $E(H)$, we can recolor $uv$ with a different available color. Since $uv$ sees at most $\deg(u)-1+1+\deg(v)-1+1\leq \Delta+d < \Delta+d+6$ colors, there must be an available color for each such edge. Thus, there exists a total-$(\Delta+d+6)$-coloring extending $H$ to $G$, contradicting that $G$ was a counterexample.
\end{proofc}

\begin{claim}\label{claim:leaf-edges}
Every vertex $v$ is adjacent to at most $d$ leaves.
\end{claim}

\begin{proofc}
Let $u$ be a vertex adjacent to a leaf. It suffices to show that $uv\in E(H)$ for any leaf $v\in N(u)$ since then $H$ having maximum degree $d$ implies the desired claim. If $v\notin V(H)$, then we can color $v$ with any color it does not see (which is at most $1$ color). Thus, we may assume $v\in V(H)$. If $u\notin V(H)$, then we can color $u$ with any color it does not see (which is at most $\Delta$ colors) and color $uv$ with any color it does not see (which is at most $\Delta+1$ colors). Thus, we may assume that $u,v\in V(H)$ and since $H$ induces a proper vertex coloring, $\varphi(v)\neq \varphi(u)$. 

Now obtain $G'$ from $G$ by removing $v$. Since $E(G')=E(G)-1$ and $|V_{[2,d+5]}|$ can only increase by at most $1$ in $G'$, we see that $q(G')< q(G)$ and by minimum counterexample we get an extension of $\varphi$ to a total-$(\Delta+d+6)$-coloring of $G'$. Since $\varphi(u)$ is precolored with some color distinct from the precoloring of $v$, we can now add back in $v$ with $\varphi(v)$ and greedily color $uv$ with any available color. Note that this is always possible since $uv$ sees at most $\Delta+1$ colors and $\Delta+d+6>\Delta+1$. Thus, we have a total-$(\Delta+d+6)$-coloring extending $H$ to $G$, contradicting that $G$ was a counterexample. 
\end{proofc}

For $i\geq 0$, let $\tilde{F}_i$ be the faces of $G$ with exactly $i$ incident vertices of degree at least $3$. 

\begin{claim}\label{F3}
$\tilde{F}_i=\emptyset$ for $i=0, 1, 2$. 
\end{claim}

\begin{proofc} Suppose not for some $i\in\{0, 1, 2\}$ . 
Suppose first that the boundary of $f$ contains
  no cycle. This means that $G$ is a forest, and $f$ is its one
  face. In particular, $G$ is bipartite. But then we may apply \cref{thm: planar-bip-prelim} and we are done. So we may assume that the boundary of $f$ contains a cycle. But then this boundary contains at least three vertices of degree at least three, since $V_2=\emptyset$ by Claim \ref{2tplus1empty}, contradiction. 
\end{proofc}

We now introduce a discharging argument. We first define an initial charge $\mu$ of $x\in V(G)\cup F(G)$. Initially let \[\mu(x) = \begin{cases}
    3\deg(x)-6 \hspace{.5cm}&\text{if $x\in V(G)$}\\
    -6 &\text{if $x\in F(G)$}\\
\end{cases}\]
By Euler's formula we see that \[\mu(G) := \sum_{v\in V(G)}\mu(v)+\sum_{f\in F(G)}\mu(f) = 6|E(G)|-6|V(G)|-6|F(G)| = -12.\]

Note leaves are the only vertices to start with a negative charge. We now present the discharging rules.

\bigskip

\begin{mdframed}
\textit{Discharging rules}
\begin{enumerate}[label=(S\arabic*)]
\item\label{S1}  For each $m$, every face $f \in \tilde{F}_m$ receives charge $\frac{6}{m}$ from each vertex of degree at least 3 on its boundary.
\item\label{S2} Every vertex $v \in V_1$ receives charge $3$ from its neighbor.
\end{enumerate}
\end{mdframed}

\medskip

Let $\mu'(x)$ be the final charge of $x\in V(G)\cup F(G)$ after \ref{S1} and \ref{S2} have been applied, and $\mu'(G)$ be the sum of all these final charges in $G$. Then $\mu'(G) = \mu(G) = -12$. We complete our proof by showing that $\mu'(G) \geq 0$, which is a contradiction.

First consider a face $f$. By Claim \ref{F3}, $f\in F_m$ for
$m\geq 3$. So by \ref{S1} (the only rule
affecting $f$), $\mu'(f)=\mu(f)+ m (\frac{6}{m})=-6 + 6 = 0$.

It remains only to consider the final charge of an arbitrary vertex $v$. If $v\in V_1$, then by \ref{S2}, we get
$\mu'(v)=\mu(v)+3=-3+3=0$. By Claim \ref{2tplus1empty}, we may now assume that $\deg_G(v)\geq d+6$.

Suppose $v$ lies on the boundary of $x$ distinct faces and is incident to $y$ leaves. We know that $x$ is no more than  $\deg_G(v)-y$,  so $x+y \leq \deg_G(v)$. We also know that $y \leq d$ by Claim \ref{claim:leaf-edges}. By doubling the first inequality and adding the result to the second inequality we get
\begin{equation}\label{2x} 2x+3y \leq 2\deg_G(v)+d.
\end{equation}
Since $\tilde{F}_i = \emptyset$ for any $0\leq i\leq 2$ by Claim \ref{F3}, each of the $x$ distinct faces incident to $v$ has at least $3$ vertices of degree at least 3 on their boundary. This means that each of these $x$ faces takes charge at most  2 from $v$, according to \ref{S1}. Each of the $y$ leaves incident to $v$ take exactly 3 from $v$, according to \ref{S2}. Hence by inequality (\ref{2x}), 
\begin{equation}\label{ab}
  \mu'(v)=3\deg_G(v)-6 -(2x+3y) \geq \deg_G(v)-6-d\geq (d+6)-d \geq 0.
\end{equation}
Thus, every vertex and face has nonnegative final charge, a contradiction which completes the proof.
\end{proof}

In the case when $H$ is a totally precolored matching ($d=1$) we get the following corollary of \cref{thm:Planar-Delta+t}

\begin{corollary}\label{D+7}
    Let $G$ be a planar graph with maximum degree $\Delta$ and let $H$ be a subgraph of $G$. If $H$ is a matching then any total-$(\Delta+7)$-coloring of $H$ in $G$ can be extended to $G$. 
\end{corollary}

\section{Proof of \cref{thm: D+3-no-distance}}\label{sec: thm1.3}

We reduce the color palette in \cref{D+7} from $(\Delta+7)$ to $(\Delta+3)$ with the following theorem.

\begin{theorem}\label{thm: No-Distance-K1-K2}
Let $3\leq t\leq 6$, let $G$ be a planar graph with maximum degree $\Delta\geq 31-t$, and let $H$ be a subgraph of $G$. If $H$ is a subgraph with maximum degree one, then any total-$(\Delta+t)$-coloring of $H$ in $G$ can be extended to $G$. 
\end{theorem}

Note that setting $t=1$ in the above statement gives a slightly strengthened version of Theorem \ref{thm: D+3-no-distance} (since here we allow $H$ to also contain isolates).

\begin{proof}
Suppose that $G$ is a counterexample for which $t$ is maximum and, subject to that, the quantity $q(G) = 3|E(G)|+|V_{[2,t-1]}|$ is minimum. Then there exists a subgraph $H$ with maximum degree one such that some total-$(\Delta+t)$-coloring $\varphi_0$ of $H$ in $G$ cannot be extended to $G$.

Consider the following claim, whose statements are very similar to claims we have already established within other proofs in this paper; the previous claims are listed in brackets below. Here, as in our previous work, we define a vertex $v$ to be \emph{high} if $\deg_G(v)\geq \lceil \frac{\Delta+t}{2} \rceil$.

\begin{claim}\label{obs: REPEATS}
The following hold.
\begin{enumerate}
\item\label{item: all but one edge} If $e\in E(G)\setminus E(H)$, then $\varphi_0$ can be extended to $G-e$. (\cref{clm: Gprime} of \cref{thm: 111t})
    \item\label{item: uv degree} Let $uv\in E(G)\setminus E(H)$ such that $u\notin V(H)$ and $u$ is low. Then $\deg(u)+\deg(v)\geq \Delta+t$. Moreover, every edge $uv\in E(G)\setminus E(H)$ has at least one high endpoint. (\cref{obs: ALL CLIQUES - deg sum} of \cref{thm: 111t})
    \item\label{item: V1} $V_1 \subseteq V(H)$ and for any $v\in V_{[t,4]}$, $N(v)\subseteq V_{[\Delta-4+t,\Delta]}$ (\cref{clm: Vt}(1) and (3) of \cref{thm: 111t})
    \item\label{item: V2-Vt-1} $V_{[2,t-1]}=\emptyset$ (\cref{2tplus1empty} of \cref{thm:Planar-Delta+t}).
\end{enumerate}
\end{claim}

To establish Claim \ref{obs: REPEATS}(1), we can use an essentially identical proof to the one indicated in the bracket. To see this, note first that if $G'$ is a proper subgraph of $G$, then $q(G')<q(G)$ (since any edges removed create at most $2$ new vertices in 
$V_{[2,t-1]}$), so we may apply minimality. When $t$ achieves its upper bound in the other proof, we applied Theorem \ref{PlanarDelta+3}, but here we can use Corollary \ref{D+7} instead. To establish Claim \ref{obs: REPEATS}(2), the proof of the bracketed claim works verbatim, except that instead of relying on ``Claim 1'' there, we instead rely on the analog we just established in Claim \ref{obs: REPEATS}(1). Similarly, the fact that $V_1\subseteq V(H)$ for Claim \ref{obs: REPEATS}(3) is verbatim from Claim \ref{clm: Vt}(1), except we replace ``Claim 2'' there by our analog in Claim \ref{obs: REPEATS}(2). The second part of Claim \ref{obs: REPEATS}(3) is only relevant when $t\in \{3,4\}$. For $t=3$ the result follows directly from Claim \ref{obs: REPEATS}(2). For $t=4$ we follow the same argument as Claim \ref{clm: Vt}(3), replacing ``Claim 1'' and ``Claim 2'' as above. Finally, to establish Claim \ref{obs: REPEATS} (4), we look to the proof of the bracketed result, where $\Delta+ d+5$ is replaced by $\Delta+t-1$ here, and we note that $t-1 \geq  \delta(H)+1$. 

Our next claim is similar to Claim 4 in the proof of Theorem \ref{thm: 111t}, but with enough differences that we write a full argument here for clarity.

\begin{claim}\label{obs: global pot}
If $t\in \{3,4\}$, then $|V_{[t,4]}|<|V_{[\Delta-4+t, \Delta]}|$
\end{claim}

\begin{proofc}
Let $B$ be the bipartite graph induced by parts $V_{[t,4]}$ and $V_{[\Delta-4+t,\Delta]}$ in $G-E(H)$. Note that by \cref{obs: REPEATS}(\ref{item: V1}) all but at most one edge incident to $v\in V_{[t,4]}$ is in $E(B)$. If $B$ is acyclic, then $2|V_{[t,4]}|\leq |E(B)|<|V(B)| = |V_{[t,4]}|+|V_{[\Delta-4+t, \Delta]}|$ as desired. 

Now suppose for contradiction that $B$ contains a cycle $C$. The edges of $C$ are all uncolored, but vertices on $C$ can be precolored. If $v\in V(C)$ and $v\not\in V(H)$, then note that $v$ sees at most $\Delta$ other colors, since no edge incident to $v$ is precolored (precolored edges have both ends precolored as well). So the precoloring $\varphi_0$ of $H$ can be extended to $H'=H\cup V(C)$. Moreover, $H'$ still has maximum degree one, since it is obtained from $H$ by adding isolates. So if we let $G'=G\setminus E(C)$, we can use \cref{obs: REPEATS}(\ref{item: all but one edge}) to get a total-$(\Delta+t)$-coloring $\varphi$ extending $\varphi_0$ to $G'$.
\color{black} It remains then only to color the edges of $C$.

Let $uv\in E(C)$, with $u\in V_{[t,4]}, v\in V_{[\Delta-4+t,\Delta]}$.
If $t=3$, the edge $uv$ sees at most $1+\Delta-2+2=\Delta+1$ colors in $\varphi$ (one edge incident to $u$, $\Delta-2$ edges incident to $v$, and the two vertices $u, v$ themselves). Since $\Delta+t=\Delta+3$ here, there are at least $2$ color choices for each such edge, so since even cycles are $2$-edge-choosable, we can extend $\varphi$ to $G$, contradiction. We can make a similar argument when $t=4$. In that case, each edge sees at most $2+\Delta-2+2=\Delta+2$ colors in $\varphi$. So, since $\Delta+t=\Delta+4$, there are at least $2$ color choices for each such edge, and the $2$-edge-choosability of $C$ means we can extend $\varphi$ to $G$, contradiction.
\end{proofc}

We say a leaf in $G$ is a \emph{high-leaf} (\emph{low-leaf}) if its neighbor is a high (low) vertex. Given a high vertex  $v\in V(G)\setminus V(H)$,  we say a face $f$ is \emph{needy with respect to $v$} if $v$ is incident to $f$ and either of the following is true:
\begin{enumerate}
    \item $f$ is a $3$-face and $v$ is its only incident high vertex, or; 
    \item $f$ is a face with at least two incident low-leaves $u,u'$ with neighbors $w,w'$, respectively, such that $w,w'$ are neighbors of $v$ on $f$.
\end{enumerate}
We say that a needy face is either \emph{type 1} or \emph{type 2} according to whether it satisfies (1) or (2), respectively See Figure \ref{fig:needyfaces}. Given a high vertex $v\in V(G)\setminus V(H)$, let $\eta(v)$ denote the number of faces that are needy with respect to $v$.

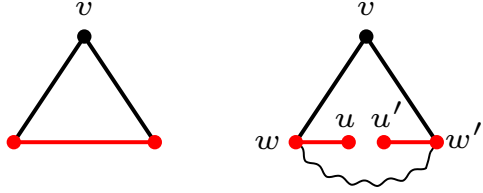
\begin{figure}[h!]
    \centering
    \begin{tikzpicture}









\begin{scope}[shift = {(5,0)}]
\draw[red, thick] (-.5,.25) to (-.125,.25);
\draw[red, thick] (.5,.25) to (.125,.25);
\draw[thick] (.5,.25) to (0,1);
\draw[thick] (-.5,.25) to (0,1);
\draw[decorate, decoration={snake, amplitude=0.2mm, segment length=1.8mm}, bend right = 60] (-.5,.25) to (.5,.25);

\draw[red, fill=red] (-.125,.25) circle (1.25pt);
\draw[red, fill=red] (.125,.25) circle (1.25pt);
\draw[fill = black] (0,1) circle (1.25pt);
\draw[red, fill=red] (-.5,.25) circle (1.25pt);
\draw[red, fill=red] (.5,.25) circle (1.25pt);

\node at (0,1.2) {\tiny{$v$}};
\node at (-.7,.25) {\tiny{$w$}};
\node at (.7,.3) {\tiny{$w'$}};
\node at (-.15,.4) {\tiny{$u$}};
\node at (.15,.45) {\tiny{$u'$}};
\end{scope}

\begin{scope}[shift = {(3,0)}]
\draw[thick] (.5,.25) to (0,1);
\draw[thick] (-.5,.25) to (0,1);
\draw[red, thick] (-.5,.25) to (.5,.25);

\draw[fill = black] (0,1) circle (1.25pt);
\draw[red, fill=red] (-.5,.25) circle (1.25pt);
\draw[red, fill=red] (.5,.25) circle (1.25pt);

\node at (0,1.2) {\tiny{$v$}};
\end{scope}

\end{tikzpicture}
    \caption{Needy faces with respect to $v$ of type 1 (left) and type 2 (right). Precolored vertices and edges indicated in red.}
    \label{fig:needyfaces}
\end{figure}

\begin{claim}\label{obs: needy}
Let $v\in V(G)\setminus V(H)$ be a high vertex. Then $\eta(v)\leq \frac{1}{2}\deg(v)$.    
\end{claim}

\begin{proofc}
Suppose for contradiction that $v\in V(G)\setminus  V(H)$ and $\eta(v)>\frac{1}{2}\deg(v)$. Then by the pigeonhole principle, there must be two adjacent needy faces with respect to $v$, say $f$ and $f'$. If $f$ and $f'$ are of type 2, this would require a neighbor of $v$ to be adjacent to two leaves. By \cref{obs: REPEATS}(\ref{item: V1}),  since $V_1\subseteq V(H)$, this would contradict $H$ being a matching. So we may assume without loss of generality that $f$ is a $3$-face where $v$ is the only high vertex. 
Let $u, v , w$ be the vertices of $f$ where $u,w$ are both low and $u$ is incident to $f'$. By \cref{obs: REPEATS}(2), $uw\in E(H)$. 

Now consider $f'$. If $f'$ is of type $2$, then $u$ must also have a neighboring leaf $z$ on $f'$ which means $uz\in E(H)$ by \cref{obs: REPEATS}(2), contradicting that $H$ is a matching. If $f'$ is of type $1$ where $f' = vuw'$ then by definition $u$ and $w'$ are both low and thus $uw'\in E(H)$ by \cref{obs: REPEATS}(2) again contradicting that $H$ is a matching.
In either case we arrive at a contradiction since $H$ is a matching.
\end{proofc}

\begin{claim}\label{obs: D+t - lengthofface}
If $f\in F(G)$ is incident with $x\geq 1$ leaves, then $\ell(f)\geq x+4$. 
\end{claim}

\begin{proofc}
We prove this by induction on $x$. If $x = 1$, then any face $f$ with one leaf $v$ must have length at least $5$ since $G$ is simple. Now suppose a face $f$ has $x\geq 2$ leaves. Let $v$ be a leaf on $f$ and let $f'$ be the face obtained by removing $v$. By \cref{obs: REPEATS}(\ref{item: V2-Vt-1}), $V_2=\emptyset$ so when we remove $v$ from $f$, we yield a face $f'$ incident with one less leaf. By induction $\ell(f')\geq x-1+4 = x+3$. By adding $v$ back in we see that $\ell(f) = \ell(f')+2 = x+5 \geq x+4$ as desired.
\end{proofc}


 We now define an initial charge $\mu$ of $x\in V(G)\cup F(G)$ and to a \emph{global pot} $\mathcal{P}$. Initially let $\mu(\mathcal{P})=0$ and let \[\mu(x) = \begin{cases}
    \deg(x)-4 \hspace{.5cm}&\text{if $x\in V(G)$}\\
    \ell(x)-4 &\text{if $x\in F(G)$}\\
\end{cases}\]
By Euler's formula we see that \[\mu(G) = \mu(\mathcal{P})+\sum_{v\in V(G)}\mu(v)+\sum_{f\in F(G)}\mu(f) = 4|E(G)|-4|V(G)|-4|F(G)| = -8.\]
We now present the discharging rules.

\medskip    

\begin{mdframed}
    
\textit{Discharging rules}

\begin{enumerate}[label=(T\arabic*)]
    \item \label{T1} Every high-leaf receives charge $1$ from its incident face and receives charge $2$ from its neighbor.
    \item \label{T2} Every low-leaf receives charge $1$ from its incident face $f$, receives charge $1$ from its neighbor and receives charge $\frac{1}{2}$ from both of its second neighbors on $f$.
    \item \label{T3} Every $3$-face receives charge $\frac{1}{y}$ from its $y$ incident high vertices (note $y\geq 1$).
    \item \label{T4} If $t\in \{3,4\}$, every vertex $v\in V_{[t,4]}$ receives charge $2$ from $\mathcal{P}$ and every vertex $v\in V_{[\Delta-4+t,\Delta]}$ gives charge $2$ to $\mathcal{P}$.
\end{enumerate}
\end{mdframed}

\medskip

Let $\mu'(x)$ and $\mu'(\mathcal{P})$ be the final charge of $x\in V(G)\cup F(G)$ and $\mathcal{P}$ after \ref{T1}-\ref{T4} have been applied. We prove via a series of three claims that $\mu'\geq 0$, contradicting $\mu'(G) = \mu(G) = -8$ and completing the proof.  Our proof consists of one claim each for $\mathcal{P}$, faces, and vertices, respectively.

\begin{claim}
$\mathcal{P}$ has nonnegative final charge.
\end{claim}

\begin{proofc}
The only rule affecting $\mathcal{P}$ is \ref{T4}. By \cref{obs: global pot},  $\mu'(\mathcal{P})\geq \mu(\mathcal{P})+2|V_{[\Delta-4+t,\Delta]}| -  2|V_{[t,4]}|>0$. 
\end{proofc}

\begin{claim}
    Every face has nonnegative final charge.
\end{claim}

\begin{proofc}
Let $f\in F(G)$. If $f\in F_3$, then $f$ cannot have any incident leaves. Thus, $f$ is only affected by \ref{T3} and $\mu'(f) = \mu(f) + y \cdot \frac{1}{y} = 0$ where $y$ is the number of high vertices incident to $f$. We reiterate that $y\geq 1$ since every $3$-face can have at most $1$ edge in $E(H)$ and hence there is an edge which is not precolored and must be incident to a high vertex. Next suppose $\ell(f)\geq 4$. If $f$ has no incident leaves, $f$ is not affected by any discharging rules and thus, $\mu'(f) = \mu(f) = \ell(f)-4\geq 0$. Finally, suppose $f$ is incident to $x\geq 1$ leaves. By \ref{T1} and \ref{T2}, $\mu'(f) = \mu(f) - x = \ell(f)-4-x$. By \cref{obs: D+t - lengthofface}, $\ell(f)\geq 4+x$ and thus, $\mu'(f)\geq 0$.
\end{proofc}

\begin{claim}
    Every vertex has nonnegative final charge. 
\end{claim}

\begin{proofc}
    Let $v\in V(G)$. First suppose $v$ is a leaf. If $v$ is a high-leaf, by \ref{T1}, $\mu'(v) = \mu(v)+1+2 = \deg(v)-4+3 = 0$. If $v$ is a low-leaf, by \ref{T2}, $\mu'(v)= \mu(v)+1+1+2\cdot\frac{1}{2} = \deg(v)+3 = 0$. 

Next suppose $v$ is low. Note that since $v$ is low, $v$ does not give any charge to $\mathcal{P}$ in \ref{T4}. By \ref{T2}, $\mu'(v) \geq \mu(v) - 1 = \deg(v)-5$ and thus $\mu'(v)\geq 0$ if $\deg(v)\geq 5$. If $2\leq \deg(v)\leq 4$ then since $V_{[2,t-1]} = \emptyset$ by \cref{obs: REPEATS}(\ref{item: V2-Vt-1}), $\deg(v)\in \{3,4\}$ and  $t\in \{3,4\}$. Thus, by \ref{T4}, $\mu'(v)\geq \mu(v)+2-1 = \deg(v)-3\geq 0$. 

Next suppose $v$ is high and $v$ is not affected by \ref{T4}. If $v\in V(H)$, then every $3$ face which $v$ is incident to, must have at least $2$ high vertices incident to it. Thus, by \ref{T1} and \ref{T3} $\mu'(v)\geq \mu(v) - (2+\frac{1}{2}(\deg(v)-1)) = \frac{1}{2}\deg(v)-\frac{11}{2}$. Thus, $\mu'(v)\geq 0$ as long as $\deg(v)\geq 11$. Since $\Delta\geq 31-t$, and $v$ is high, $\deg(v)\geq \lceil \frac{\Delta+t}{2}\ \rceil \geq 16$. If $v\notin V(H)$, then by \ref{T2} and \ref{T3}, $v$ gives a total charge of $1$ for any needy fact with respect to $v$ (gives 1 to a needy face of type 1 and $\frac{1}{2}$ to both leaves in a needy face of type 2). For any other incident face, $v$ gives charge at most $\frac{1}{2}$. Thus, by \cref{obs: needy}, we see that $\mu'(v)\geq \mu(v) - (1\cdot \eta(v)+ \frac{1}{2}\cdot (\deg(v)-\eta(v)))\geq \frac{1}{4}\deg(v)-4$. Hence $\mu'(v)\geq 0$ so long as $\deg(v)\geq 16$. Since $v$ is high and $\Delta\geq 31-t$, $\lceil \frac{\Delta+t}{2}\rceil \geq 16$.

Finally suppose $v$ is high and is affected by \ref{T4}. Then $v$ loses an additional charge of $2$ to $\mathcal{P}$. In this case, $t\in \{3,4\}$ and $\deg(v)\geq \Delta-1$. If $v\in V(H)$, then by the same computation as the previous paragraph we see that $\mu'(v)\geq \frac{1}{2}\deg(v)-\frac{11}{2}-2\geq \frac{1}{2}\Delta-8\geq 0$ since $\Delta\geq 25$. Similarly if $v\notin V(H)$, using the same calculation in the previous paragraph, $\mu'
(v) \geq \frac{1}{4}\deg(v)-6\geq \frac{1}{4}\Delta -\frac{25}{4}\geq 0$ since $\Delta\geq 25$.
\end{proofc}

\end{proof}

\bibliographystyle{abbrv}

\end{document}